\documentclass[10pt]{amsart}
\usepackage{amssymb}
\usepackage{amscd}

\theoremstyle{definition}
\newtheorem{thm}{Theorem}[section]
\newtheorem{lem}[thm]{Lemma}
\newtheorem{prp}[thm]{Proposition}
\newtheorem{dfn}[thm]{Definition}
\newtheorem{cor}[thm]{Corollary}

\newtheorem{rmk}[thm]{Remark}
\newtheorem{ntn}[thm]{Notation}
\newtheorem{exa}[thm]{Example}

\newcommand{\beq}{\begin{equation}}
\newcommand{\eeq}{\end{equation}}
\newcommand{\beqr}{\begin{eqnarray*}}
\newcommand{\eeqr}{\end{eqnarray*}}
\newcommand{\bal}{\begin{align*}}
\newcommand{\eal}{\end{align*}}
\newcommand{\bei}{\begin{itemize}}
\newcommand{\eei}{\end{itemize}}

\newcommand{\af}{\alpha}
\newcommand{\bt}{\beta}
\newcommand{\gm}{\gamma}
\newcommand{\dt}{\delta}
\newcommand{\ep}{\varepsilon}
\newcommand{\zt}{\zeta}
\newcommand{\et}{\eta}

\newcommand{\te}{\theta}
\newcommand{\ld}{\lambda}
\newcommand{\sm}{\sigma}

\newcommand{\ph}{\varphi}
\newcommand{\ps}{\psi}
\newcommand{\rh}{\rho}
\newcommand{\om}{\omega}
\newcommand{\ta}{\tau}

\newcommand{\Gm}{\Gamma}

\newcommand{\Ld}{\Lambda}

\newcommand{\Q}{{\mathbb{Q}}}
\newcommand{\Z}{{\mathbb{Z}}}
\newcommand{\R}{{\mathbb{R}}}
\newcommand{\C}{{\mathbb{C}}}
\newcommand{\N}{{\mathbb{N}}}

\newcommand{\id}{{\mathrm{id}}}
\newcommand{\ev}{{\mathrm{ev}}}

\newcommand{\dist}{{\mathrm{dist}}}

\newcommand{\diag}{{\mathrm{diag}}}
\newcommand{\supp}{{\mathrm{supp}}}

\newcommand{\spn}{{\mathrm{span}}}

\newcommand{\Aut}{{\mathrm{Aut}}}
\newcommand{\Ad}{{\mathrm{Ad}}}

\newcommand{\Ker}{{\mathrm{Ker}}}

\newcommand{\Zq}[1]{{\Z_{#1}}}

\newcommand{\Zqn}{\Zq{n}}
\newcommand{\Cs}[3]{C^* (\Zq{#1}, #2, #3)}

\newcommand{\CZnAa}{\Cs{n}{A}{\af}}

\newcommand{\andeqn}{\,\,\,\,\,\, {\mbox{and}} \,\,\,\,\,\,}
\newcommand{\QED}{\rule{0.4em}{2ex}}
\newcommand{\ts}[1]{{\textstyle{#1}}}
\newcommand{\ds}[1]{{\displaystyle{#1}}}
\newcommand{\ssum}[2]{{\ts{ {\ds{\sum}}_{#1}^{#2} }}}

\newcommand{\ca}{C*-algebra}

\newcommand{\ct}{continuous}
\newcommand{\pj}{projection}

\newcommand{\nbhd}{neighborhood}

\newcommand{\hm}{homomorphism}
\newcommand{\wolog}{without loss of generality}
\newcommand{\Wolog}{Without loss of generality}

\newcommand{\ifo}{if and only if}

\newcommand{\mops}{mutually orthogonal \pj s}

\newcommand{\cfn}{continuous function}
\newcommand{\hsa}{hereditary subalgebra}

\newcommand{\mvnt}{Murray-von Neumann equivalent}

\newcommand{\tRp}{tracial Rokhlin property}

\newcommand{\fd}{finite dimensional}
\newcommand{\uct}{Universal Coefficient Theorem}
\newcommand{\suct}{satisfies the \uct}

\newcommand{\rsz}[1]{\raisebox{0ex}[0.8ex][0.8ex]{$#1$}}

\renewcommand{\S}{\subset}

\newcommand{\SM}{\setminus}

\title[Simple higher dimensional noncommutative torus]{Every simple higher dimensional noncommutative torus is an AT~algebra}

\author{N.~Christopher Phillips}

\date{17~Sept.\  2006}  %

\address{Department of Mathematics, University  of Oregon,
       Eugene OR 97403-1222, USA.}

\email[]{ncp@darkwing.uoregon.edu}

\thanks{Research partially supported by
  NSF grants DMS~0070776 and DMS~0302401.}

\begin{document}

\setcounter{section}{-1}

\begin{abstract}
We prove that 
every simple higher dimensional noncommutative torus is an AT~algebra.
\end{abstract}

\maketitle

\section{Introduction}\label{Sec:Intro}

\indent
A higher dimensional noncommutative torus is the universal \ca\  %
generated by unitaries which commute up to specified scalars.
Thus, it is a generalization
of the rotation algebra $A_{\te}$ to more generators.
The commutation relations are determined by a real skew symmetric
matrix; see Notation~\ref{StdNtn} for a precise formulation.
In this paper, we prove that
every simple higher dimensional noncommutative torus
is an AT~algebra, that is, a direct limit of finite direct sums
of \ca s of the form $C (S^1, M_n)$ for varying values of $n.$

The first result in this direction is the
Elliott-Evans Theorem~\cite{EE} for the ordinary irrational
rotation algebras, which we use here as the initial step of an
induction argument.
Without giving a complete list of later work, we mention four
highlights.
All simple three dimensional noncommutative tori were shown to be
AT~algebras by Q.~Lin~\cite{LQ1}.
In arbitrary dimension,
``most'' simple higher dimensional noncommutative tori were
shown to be AT~algebras by Boca~\cite{Bc}.
Kishimoto (Corollary~6.6 of~\cite{Ks4})
obtained this result in all cases in which,
in the skew symmetric matrix giving the commutation relations,
the entries above the diagonal are rationally independent,
as well as some others.
Theorem~3.14 of~\cite{Ln14} shows that the crossed product of
$(S^1)^{d}$ by a minimal rotation is an AT~algebra; in this case,
most of the entries of the relevant skew symmetric matrix are zero.

Our proof is by induction on the number of generators.
Every higher dimensional noncommutative torus can be written
as an iterated crossed product by $\Z,$
and the proof of Kishimoto's result uses an inductive
argument which works whenever the intermediate crossed
products are all simple.
One has some choice here: different choices of the commutation
relations may well give the same \ca.
As a very simple example, one might simply write the generators
in a different order.
Unfortunately, it seems not to be possible in general to choose
commutation relations to give the same algebra, or even a
Morita equivalent algebra (see~\cite{RS}), and in such a way
that Kishimoto's method applies, or even
in such a way as to get a tensor product of algebras to which
this method applies.
However, if one allows one more kind of modification, namely the
replacements of unitary generators by integer powers of themselves,
then it is always possible to replace a noncommutative torus by a
tensor product of algebras covered by Kishimoto's method.
The new algebra isn't isomorphic, or even Morita equivalent, to
the original.
But if one replaces only one generator, the new algebra is the fixed
point algebra of a tracially approximately inner action~\cite{PhtRp1a}
of a finite cyclic group which has the \tRp~\cite{PhtRp1a}.
As proved in~\cite{PhtRp1a},
this operation thus preserves tracial rank zero.
Because $K_0$ and $K_1$ are torsion free, H.~Lin's classification
theorem for simple nuclear \ca s with tracial rank zero,
Theorem~5.2 of~\cite{Ln15}, shows
that this operation preserves the property of being an AT~algebra.
Indeed, the \tRp\  was introduced specifically for this purpose.
It is not possible to substitute the Rokhlin property
of~\cite{Iz} and~\cite{Iz2} in this argument.
As we show in Corollary~\ref{NoRP},
the relevant action does not have the Rokhlin property.

This paper is organized as follows.
Section~\ref{Sec:NCT} contains various preliminaries.
There is little that is really new,
but our presentation gives the material in a convenient form
and establishes notation for the rest of the paper.
In Section~\ref{Sec:ARokhOnNCT}, we prove that if $A$ is a simple
higher dimensional noncommutative torus,
then the automorphism which multiplies one
of the standard unitary generators by $\exp (2 \pi i / n)$
generates an action  of $\Zqn$ with the \tRp.
In Section~\ref{Sec:NCTInd}, we use this result, Kishimoto's result,
and H.~Lin's classification theorem~\cite{Ln15},
to construct an inductive proof that every simple noncommutative torus
is an AT~algebra,
and we obtain several corollaries.

This paper contains the material of Sections~5 through~7
of the unpublished long preprint~\cite{PhW}.

We use the notation $\Zq{n}$ for $\Z / n \Z$;
the $p$-adic integers will not appear in this paper.
If $A$ is a \ca\  and $\af \colon A \to A$ is an automorphism
such that $\af^n = \id_A,$ then we write
$\CZnAa$ for the crossed product of
$A$ by the action of $\Zqn$ generated by $\af.$
We write $p \precsim q$ to mean that the \pj\  $p$ is \mvnt\  to a
sub\pj\  of $q,$ and $p \sim q$ to mean that
$p$ is \mvnt\  to $q.$
Also, $[a, b]$ denotes the additive commutator $a b - b a.$

I am grateful to Marc Rieffel for discussions concerning
higher dimensional noncommutative tori,
and in particular for pointing out that it was not known
whether a higher dimensional noncommutative torus is isomorphic to
its opposite algebra.
See Corollary~\ref{IsomOpp};
this remains open in the nonsimple case.
I would like to thank Hiroyuki Osaka for carefully
reading earlier versions of this paper, and catching a number
of misprints and suggesting many improvements,
and Hanfeng Li for useful comments,
including the improvement of Lemma~\ref{FPOfMult}.

\section{Higher dimensional noncommutative tori}\label{Sec:NCT}

\indent
In this section we present, in a form convenient for our purposes,
some mostly standard facts
about higher dimensional noncommutative tori
and about irrational rotation algebras.

\begin{ntn}\label{StdNtn}
Let $\te$ be a skew symmetric real $d \times d$ matrix.
The {\emph{noncommutative torus}} $A_{\te}$ is by definition~\cite{Rf2}
the universal \ca\  generated by unitaries $u_1, u_2, \dots, u_d$
subject to the relations
\[
u_k u_j = \exp (2 \pi i \te_{j, k} ) u_j u_k
\]
for $1 \leq j, k \leq d.$
(Of course, if all $\te_{j, k}$ are integers, it is not really
noncommutative.)
\end{ntn}

Some authors use $\te_{k, j}$ in the commutation relation instead.
See for example~\cite{Ks2}.

\begin{rmk}\label{CoordFree}
We note (see the beginning of Section~4 of~\cite{Rf1} and
the introduction to~\cite{RS}) that $A_{\te}$ is the universal
\ca\  generated by unitaries $u_x,$ for $x \in \Z^d,$
subject to the relations
\[
u_y u_x = \exp (\pi i \langle x, \te (y) \rangle ) u_{x + y}
\]
for $x, \, y \in \Z^d.$

It follows that if $B \in {\mathrm{GL}}_d (\Z),$ and if
$B^{\mathrm{t}}$ denotes the transpose of $B,$ then
$A_{B^{\mathrm{t}} \te B} \cong A_{\te}.$
That is, $A_{\te}$ is unchanged if $\te$ is rewritten in terms of
some other basis of $\Z^d.$
\end{rmk}

\begin{rmk}\label{AlgFromBichar}
Let $\af$ be a skew symmetric real bicharacter on $\Z^d,$ that is,
a $\Z$-bilinear function $\af \colon \Z^d \times \Z^d \to \R$
such that $\af (x, y) = - \af (y, x)$ for all $x, \, y \in \Z^d.$
For any basis $(b_1, b_2, \dots, b_d)$ of $\Z^d,$
there is a unique skew symmetric real $d \times d$ matrix $\te$
such that
\[
\af \left( \ssum{k = 1}{d} x_k b_k, \, \ssum{k = 1}{d} y_k b_k \right)
  = \ssum{j, k = 1}{d}  x_j \te_{j, k} y_k
\]
for all $x, \, y \in \Z^d.$
We define $A_{\af} = A_{\te}.$
Remark~\ref{CoordFree} shows that this \ca\  is independent of the
choice of basis.
\end{rmk}

\begin{rmk}\label{Restrict}
Let $\af$ be a skew symmetric real bicharacter on $\Z^d,$ and let
$H \S \Z^d$ be a subgroup.
Then $H \cong \Z^m$ for some $m \leq d.$
By abuse of notation, we write $\af |_H$ for the restriction of
$\af$ to $H \times H \S \Z^d \times \Z^d.$
There is a noncommutative torus $A_{\af |_H}$ by
Remark~\ref{AlgFromBichar}, which does not depend on the choice
of the isomorphism $H \cong \Z^m.$

For a skew symmetric real $d \times d$ matrix $\te$ and a subgroup
$H \S \Z^d$ with a specified ordered basis, we write
$\te |_H$ for the matrix in that basis of the restriction to $H$ of
the real bicharacter $(x, y) \mapsto \langle x, \te y \rangle.$
For subgroups such as $\Z^m \times \{ 0 \}$ or
$\Z^m \times \{ 0 \} \times \Z^l,$ we use without comment the
obvious basis.
\end{rmk}

We formalize a remark made in 1.7 of~\cite{El0},
according to which every
noncommutative torus can be obtained as a repeated crossed
product by $\Z.$

\begin{lem}\label{ItCrPrd}
Let $\af$ be a skew symmetric real bicharacter on $\Z^d.$
Then there is an automorphism
$\ph$ of $A_{\af |_{\Z^{d - 1} \times \{ 0 \} }}$ which is homotopic
to the identity and such that
\[
A_{\af}
 \cong C^* ( \Z, \, A_{\af |_{\Z^{d - 1} \times \{ 0 \} }}, \, \ph).
\]
\end{lem}

\begin{proof}
Let $\te$ be the matrix of $\af$ in the standard basis.
Let $\bt = \af |_{\Z^{d - 1} \times \{ 0 \} }.$
Then the matrix of $\bt$ is
$( \te_{j, k} )_{1 \leq j, k \leq d - 1}.$
Let $u_1, \, u_2, \, \dots, \, u_{d - 1}$ be the standard
generators of $A_{\bt}.$
Then $\ph$ is determined by
$\ph (u_j) = \exp (2 \pi i \af_{j, d}) u_j.$
It is clear that $\ph$ is homotopic to the identity.
\end{proof}

The following definition is essentially from Section~1.1 of~\cite{Sl}.

\begin{dfn}\label{NondegDfn}
The skew symmetric real $d \times d$ matrix $\te$ is
{\emph{nondegenerate}} if whenever $x \in \Z^d$ satisfies
$\exp (2 \pi i \langle x, \te y \rangle ) = 1$ for all $y \in \Z^d,$
then $x = 0.$
Otherwise, $\te$ is {\emph{degenerate}}.
We similarly refer to degeneracy and nondegeneracy of
a skew symmetric real bicharacter on $\Z^d.$
\end{dfn}

\begin{lem}\label{NondegCond}
Let $\te$ be a skew symmetric real $d \times d$ matrix.
Then $\te$ is degenerate \ifo\  there
exists $x \in \Q^d \SM \{ 0 \}$ such that
$\langle x, \, \te y \rangle \in \Q$ for all $y \in \Q^d.$
\end{lem}

\begin{proof}
If $\te$ is degenerate, choose $w \neq 0$ such that
$\exp (2 \pi i \langle w, \te y \rangle ) = 1$ for all $y \in \Z^d.$
Then $\langle w, \te y \rangle \in \Z$ for all $y \in \Z^d.$
If now $y \in \Q^d$ is arbitrary,
then there exists $m \in \Z \SM \{ 0 \}$ such that $m y \in \Z^d.$
So
\[
\langle w, \te y \rangle
 = \ts{\frac{1}{m}} \langle w, \te (m y) \rangle
 \in \ts{\frac{1}{m}} \Z \S \Q.
\]

Conversely, assume $x \in \Q^d \SM \{ 0 \}$ and
$\langle x, \, \te y \rangle \in \Q$ for all $y \in \Q^d.$
Choose $m \in \Z$ with $m > 0$ such that
$m \langle x, \, \te \dt_k \rangle \in \Z$ for $1 \leq k \leq d.$
Then $m x \neq 0$ and
$\exp (2 \pi i \langle m x, \te y \rangle ) = 1$ for all $y \in \Z^d.$
\end{proof}

\begin{lem}\label{ConjByRat}
Let $\te$ be a skew symmetric real $d \times d$ matrix.
Let $B \in {\mathrm{GL}}_d (\Q).$
Then $B^{\mathrm{t}} \te B$ is nondegenerate
\ifo\  $\te$ is nondegenerate.
\end{lem}

\begin{proof}
It suffices to prove one direction.
Suppose $\te$ is degenerate.
By Lemma~\ref{NondegCond}, there is
$x \in \Q^d \SM \{ 0 \}$ such that
$\langle x, \, \te y \rangle \in \Q$ for all $y \in \Q^d.$
Then $B^{-1} x \in \Q^d \SM \{ 0 \}$ and
\[
\langle B^{-1} x, \, B^{\mathrm{t}} \te B y \rangle
  = \langle x, \, \te B y \rangle \in \Q
\]
for all $y \in \Q^d.$
So $B^{\mathrm{t}} \te B$ is degenerate.
\end{proof}

The following result is well known.

\begin{thm}\label{Simplicity}
The \ca\  $A_{\te}$ of Notation~\ref{StdNtn} is simple \ifo\  $\te$
is nondegenerate.
Moreover, if $A_{\te}$ is simple it has a unique tracial state.
\end{thm}

\begin{proof}
If $\te$ is nondegenerate, then $A_{\te}$ is simple by Theorem~3.7
of~\cite{Sl}.
(Note the standing assumption of nondegeneracy throughout
Section~3 of~\cite{Sl}.)

When $A_{\te}$ is simple, the proof of Lemma~3.1 of~\cite{Sl}
shows that $A_{\te}$ can have at most one tracial state.
Existence of a tracial state is well known, or can be
obtained from Lemma~\ref{ItCrPrd} by induction on $n.$

If $\te$ is degenerate, then we follow 1.8 of~\cite{El0}.
Choose $n \in \Z^d \SM \{ 0 \}$ such that
$\exp (2 \pi i \langle n, \te y \rangle ) = 1$ for all $y \in \Z^d.$
Then $v = u_1^{n_1} u_2^{n_2} \cdots u_d^{n_d}$
is a nontrivial element of the center of $A_{\te},$
which is therefore not simple.
\end{proof}

Finally, for reference and to establish notation,
we recall several facts about the ordinary rotation algebras
(the case $d = 2$).
(As far as we know, the last lemma has not appeared before,
but its proof is easy.)
We will consider various embeddings of rotation algebras
into higher dimensional noncommutative tori.
Therefore, for $\et \in \R$ we let $v_{\et}$ and $w_{\et}$ denote
the standard unitary generators of $A_{\et},$
satisfying $w_{\et} v_{\et} = \exp (2 \pi i \et) v_{\et}w_{\et}.$
In this way,
we avoid confusion with the generators $u_1, u_2, \dots, u_d$
of a higher dimensional noncommutative torus $A_{\te}.$

The proof of the next theorem
is contained in Theorem~1.1 and Proposition~1.3 of~\cite{AP}.
Also see Corollary~3.6 and Definition~3.3 of~\cite{RfF}.
We refer to~\cite{Dx} for information on \ct\  fields of \ca s.
See especially Sections~10.1 and~10.3.

\begin{thm}\label{FieldOfRotAlg}
For $\et \in \R$ let $A_{\et}$ be the rotation algebra,
with generators as described above.
Let $A$ be the \ca\  of the discrete Heisenberg group,
which is the universal \ca\  generated by unitaries $v, w, z$
subject to the relations
\[
w v = z v w, \,\,\,\,\,\, z v = v z, \andeqn z w = w z.
\]
Then there is a \ct\  field of \ca s over $S^1$ whose fiber
over $\exp (2 \pi i \et)$ is $A_{\et},$ whose \ca\  of \ct\  sections
is $A,$ and such that the
evaluation map $\ev_{\et} \colon A \to A_{\et}$ of sections
at $\exp (2 \pi i \et)$ is determined by
\[
\ev_{\et} (v) = v_{\et}, \,\,\,\,\,\, \ev_{\et} (w) = w_{\et},
\andeqn \ev_{\et} (z) = \exp (2 \pi i \et) \cdot 1.
\]
\end{thm}

Since we will only formally deal with one \ct\  field in this paper,
the following notation is unambiguous.

\begin{ntn}\label{SectionNot}
For a subset $E \S S^1,$ we let $\Gm (E)$ be the set of \ct\  sections
of the \ct\  field of Theorem~\ref{FieldOfRotAlg} over $E.$
(See 10.1.6 of~\cite{Dx}.)

For any such section $a,$ we further write
$a (\et)$ for $a (\exp (2 \pi i \et)).$
No confusion should arise.
\end{ntn}

\begin{lem}\label{ContOfTrace}
Let the notation be as in Theorem~\ref{FieldOfRotAlg}
and Notation~\ref{SectionNot}.
Let $\ta_{\et}$ be the standard trace on $A_{\et},$ satisfying
$\ta (1) = 1$ and
$\ta_{\et} (v_{\et}^m w_{\et}^n) = 0$ unless $m = n = 0.$
(If $A_{\et}$ is viewed as a crossed product by rotation
on the circle,
then $\ta_{\et}$ comes from normalized Haar measure on the circle.)
Let $U \S S^1$ be an open set, and let $a \in \Gm (U).$
Then $\et \mapsto \ta_{\et} (a (\et))$ is \ct.
\end{lem}

\begin{proof}
We check continuity at $\et_0.$
Choose a \cfn\  $h \colon S^1 \to [0, 1]$ such that $\supp (h) \S U$
and such that $h = 1$ on a \nbhd\  of $\et_0.$
Then it suffices to consider the \ct\  section $h a$ in place of $a.$
Now $h a$ is the restriction to $U$ of a \ct\  section $b$ defined on
all of $S^1,$ satisfying $b (\zt) = 0$ for $\zt \not\in U.$
Accordingly, we may restrict to the case $U = S^1.$
Then $a \in A.$

{}From the formulas
\[
\ev_{\et} (v) = v_{\et}, \,\,\,\,\,\, \ev_{\et} (w) = w_{\et},
\andeqn \ev_{\et} (z) = \exp (2 \pi i \et) \cdot 1
\]
and the definition of $\ta_{\et},$
it is immediate that if $b$ is any (noncommutative) monomial in
$v,$ $w,$ $z,$ and their adjoints,
then $\et \mapsto \ta_{\et} (b (\et))$ is \ct.
Therefore the same holds for any noncommutative polynomial,
and hence for any norm limit of noncommutative polynomials,
including $a.$
\end{proof}

\begin{lem}\label{2D}
Let the notation be as in Theorem~\ref{FieldOfRotAlg}
and Lemma~\ref{ContOfTrace}.
Let $\et \in \R \SM \Q.$
Let $n \in \N,$ let $\om = \exp (2 \pi i / n),$ and let
$\af \colon A_{\et} \to A_{\et}$ be the unique automorphism
satisfying $\af (v_{\et}) = \om v_{\et}$ and $\af (w_{\et}) = w_{\et}.$
Then for every $\ep > 0$ there exist
\mops\  $e_0, e_1, \dots, e_{n - 1}$ such that (with $e_n = e_0$)
we have $\af (e_j) = e_{j + 1}$ for $0 \leq j \leq n - 1,$ and such that
$1 - n \ta_{\et} (e_0) < \ep.$
\end{lem}

\begin{proof}
Set $\ep_0 = \frac{1}{4 n} \ep.$
Let $f \colon S^1 \to [0, 1]$ be a \cfn\  such that
$\supp (f)$ is contained in the open arc from $1$ to $\om,$
and such that $f (\zt) = 1$ for all $\zt$ in the closed arc
from $\exp (2 \pi i \ep_0)$ to
$\exp \left( 2 \pi i \left[ \frac{1}{n} - \ep_0 \right] \right).$
Then $f (v_{\et})$ is a positive element of $A_{\et}$ with
$\| f (v_{\et}) \| \leq 1$ and
$\ta_{\et} (f (v_{\et})) \geq \frac{1}{n} - 2 \ep_0.$
Since $A_{\et}$ has real rank zero (see Remark~6 of~\cite{EE},
or Theorem~1.5 of~\cite{BKR}),
there is a \pj\  $e_0$ in the \hsa\  $B$ of $A_{\et}$ generated by
$f (v_{\et})$ such that $\| e_0 f (v_{\et}) - f (v_{\et}) \| < \ep_0.$
Therefore
$\| e_0 f (v_{\et}) e_0 - f (v_{\et}) \| < 2 \ep_0.$
Since $e_0 f (v_{\et}) e_0 \leq e_0,$ it follows that
\[
\ta_{\et} (e_0) \geq \ta_{\et} (e_0 f (v_{\et}) e_0)
   > \ta_{\et} (f (v_{\et}) ) - 2 \ep_0
   \geq \ts{ \frac{1}{n} } - 4 \ep_0.
\]

We have $\af^k (f (v_{\et})) \af^l (f (v_{\et})) = 0$ for
$0 \leq k, \, l \leq n - 1$ and $k \neq l.$
Therefore $\af^k (B) \af^l (B) = \{ 0 \}$ for such $k$ and $l,$
whence also $\af^k (e_0) \af^l (e_0) = 0.$
Define $e_k = \af^k (e_0)$ for $0 \leq k \leq n - 1.$
Then $e_0, e_1, \dots, e_{n - 1}$ are \mops\  such that
$\af (e_j) = e_{j + 1}$ for $0 \leq j \leq n - 1.$
Moreover,
\[
1 - n \ta_{\et} (e_0)
  < 1 - n \left( \ts{ \frac{1}{n} } - 4 \ep_0 \right)
  = 4 n \ep_0 = \ep,
\]
as desired.
\end{proof}

\section{The tracial Rokhlin property and
  higher dimensional noncommutative tori}\label{Sec:ARokhOnNCT}

\indent
In this section, we prove that if $\te$ is nondegenerate, then the
action of $\Zqn$ which multiplies one of the standard generators
of $A_{\te}$ by a primitive $n$-th root of~$1$ has the \tRp.
We note for comparison the related result in Section~6 of~\cite{Ks2}:
if $\af \in \Aut (A_{\te})$ is of the form
$\af (u_j) = \ld_j u_j,$ with $\ld_1, \ld_2, \dots, \ld_n \in S^1,$
and if all positive powers of $\af$ are outer,
then $\af$ has the Rokhlin property.
However, our action of $\Zqn$ does not have the Rokhlin property.
See Corollary~\ref{NoRP} below.

One might hope to prove the \tRp\  fairly directly,
using one of the criteria in Section~5 of~\cite{ELP}
or Section~1 of~\cite{PhtRp1b}.
This approach does not work,
because these criteria all assume that the \ca\  involved
is already known to have tracial rank zero,
while our argument requires knowing that
the action on a particular higher dimensional noncommutative torus
has the \tRp\  before we know that it has tracial rank zero.
Instead, as is done in the proof in~\cite{BDR}
that $A_{\te}$ has real rank zero,
and analogously to Section~6 of~\cite{Ks2},
we will reduce to a construction in the ordinary
irrational rotation algebras.
The idea is to find an approximately central copy
of an ordinary irrational rotation algebra $A_{\et},$
such that the restriction to it of our action
is the one in Lemma~\ref{2D}.
Since the \pj s in $A_{\et}$ must be chosen ahead of time,
at least approximately, we must require
that $\et$ be arbitrarily close to some fixed $\et_0.$
Nondegeneracy enters through Lemma~\ref{DensityCond} below.
To obtain the correct restricted action,
we use the condition~(\ref{MakeSubalgStep2:3})
in Lemma~\ref{MakeSubalgStep2} below.
{}From then on, we roughly follow the argument used in~\cite{BKR}
to prove approximate divisibility.
We vary the arrangement slightly to make part of the argument
easily available for use elsewhere.

We recall the \tRp\  for actions of finite cyclic groups
on infinite dimensional finite simple unital \ca s.
The following result is Lemma~1.16 of~\cite{PhtRp1a},
for the case that the group is cyclic.

\begin{prp}\label{TRPCond}
Let $A$ be an infinite dimensional finite simple unital \ca,
and let $\af \in \Aut (A)$ satisfy $\af^n = \id_A.$
The action of $\Zqn$ generated by $\af$ has the
tracial Rokhlin property \ifo\  for every finite set
$F \S A,$ every $\ep > 0,$
and every nonzero positive element $x \in A,$
there are \mops\  $e_0, e_1, \dots, e_{n - 1} \in A$ such that:
\begin{enumerate}
\item\label{TRPCond:1} %
$\| \af (e_j) - e_{j + 1} \| < \ep$ for $0 \leq j \leq n - 1,$
where by convention we take the indices mod $n,$
that is, $e_n = e_0.$
\item\label{TRPCond:2} %
$\| e_j a - a e_j \| < \ep$ for $0 \leq j \leq n - 1$ and all $a \in F.$
\item\label{TRPCond:3} %
With $e = \sum_{j = 0}^{n - 1} e_j,$ the \pj\  $1 - e$ is \mvnt\  to a
\pj\  in the \hsa\  of $A$ generated by $x.$
\end{enumerate}
\end{prp}

\begin{dfn}\label{HomToInnD}
Let $\te$ be a skew symmetric real $d \times d$ matrix.
Let
\[
n = (n_1, n_2, \dots, n_d) \in \Z^d \andeqn
v = u_1^{n_1} u_2^{n_2} \cdots u_d^{n_d} \in A_{\te}.
\]
We write $\gm_n$ for the inner automorphism $\Ad (v)$ of the
noncommutative torus $A_{\te}.$
We further define a \hm\  $\sm \colon \Z^d \to (S^1)^d$ by
the formula $\sm (n)_j = \exp (2 \pi i (\te n)_j)$ for
$1 \leq j \leq d.$
(Here the expression $\te n$
is the usual action of a $d \times d$ matrix on an element of $\R^d.$)
\end{dfn}

\begin{lem}\label{HomToInnW}
Let $\te$ be a skew symmetric real $d \times d$ matrix.
With $\gm$ and $\sm$ as in Definition~\ref{HomToInnD},
we have $\gm_n (u_j) = \sm (n)_j u_j$ for $n \in \Z^d$ and
$1 \leq j \leq d.$
Moreover, if $m \in \Z^d,$ then
\[
\gm_n ( u_1^{m_1} u_2^{m_2} \cdots u_d^{m_d} )
  = \exp ( 2 \pi i \langle m, \te n \rangle )
               u_1^{m_1} u_2^{m_2} \cdots u_d^{m_d}
\]
for all $n \in \Z^d.$
\end{lem}

\begin{proof}
The first formula is the special case of the second obtained
by setting $m = \dt_j,$ the $j$-th standard basis vector of $\Z^d.$
By linearity, both formulas will follow if we check the first when
$m = \dt_j$ and $n = \dt_k.$
Since $(\te \dt_k)_j = \te_{j, k},$ this is just
the commutation relation
\[
u_k u_j u_k^* = \exp (2 \pi i \te_{j, k} ) u_j,
\]
which is the same as the one in from Notation~\ref{StdNtn}.
\end{proof}

\begin{lem}\label{DensityCond}
Let $\te$ be a skew symmetric real $d \times d$ matrix.
The \hm\  $\sm \colon \Z^d \to (S^1)^d$ of Definition~\ref{HomToInnD}
has dense range \ifo\  $\te$ is nondegenerate.
\end{lem}

\begin{proof}
Assume $\sm$ does not have dense range.
Let $H = {\overline{\sm (\Z^d) }},$ which is a proper closed subgroup
of $(S^1)^d.$
Choose a nontrivial character $\mu \colon (S^1)^d \to S^1$ whose
kernel contains $H.$
By the identification of the dual group of $(S^1)^d,$
there is $r \in \Z^d \SM \{ 0 \}$ such that
\[
\mu (\zt_1, \zt_2, \dots, \zt_d)
   = \zt_1^{r_1} \zt_2^{r_2} \cdots \zt_d^{r_d}
\]
for all $\zt \in (S^1)^d.$
Because $H \S \Ker (\mu),$ for all $n \in \Z^d$ we have
\begin{align*}
1 & = \mu ( \sm (n))
    = \exp ( 2 \pi i (\te n)_1)^{r_1}
      \exp ( 2 \pi i (\te n)_2)^{r_2} \cdots
      \exp ( 2 \pi i (\te n)_d)^{r_d}   \\
  & = \exp ( 2 \pi i \langle r, \te n \rangle).
\end{align*}
Thus $\te$ is degenerate.

Now suppose that $\te$ is degenerate.
Then we may choose $r \in \Z^d \SM \{ 0 \}$ such that
$\exp ( 2 \pi i \langle r, \te n \rangle) = 1$ for all $n \in \Z^d.$
Reversing the above calculation,
we find that the nontrivial character
\[
\mu (\zt_1, \zt_2, \dots, \zt_d)
   = \zt_1^{r_1} \zt_2^{r_2} \cdots \zt_d^{r_d}
\]
satisfies $\mu (\sm (n)) = 1$ for all $n \in \Z^d.$
Therefore $\sm$ does not have dense range.
\end{proof}

\begin{cor}\label{DensityFinInd}
Let $\te$ be a nondegenerate skew symmetric real $d \times d$ matrix.
Let $G \S \Z^d$ be a subgroup with finite index.
Let $\sm \colon \Z^d \to (S^1)^d$
be the \hm\   of Definition~\ref{HomToInnD}.
Then $\sm (G)$ is dense in $(S^1)^d.$
\end{cor}

\begin{proof}
Let $H = {\overline{\sm (G) }}.$
Let $S$ be a set of coset representatives for $G$ in $\Z^d.$
Then the sets $\sm (m) H,$ for $m \in S,$ are closed
and are pairwise equal or disjoint.
By Lemma~\ref{DensityCond}, their union is $(S^1)^d.$
Since there are finitely many of them, and since $(S^1)^d$ is
connected, it follows that all are equal to $(S^1)^d.$
\end{proof}

\begin{cor}\label{GaugeAppInn}
Let $\te$ be a nondegenerate skew symmetric real $d \times d$ matrix.
Let $\zt_1, \zt_2, \dots, \zt_d \in S^1.$
Let $\af \in A_{\te}$ be the automorphism determined by
$\af (u_j) = \zt_j u_j$ for $1 \leq j \leq d.$
Then $\af$ is approximately inner.
\end{cor}

\begin{proof}
It suffices to find, for all $\ep > 0,$
a unitary $v \in A_{\te}$ such that
$\| \af (u_j) - v u_j v^* \| < \ep$ for $1 \leq j \leq d.$
Choose $\dt > 0$ small enough that if
$(\om_1, \om_2, \dots, \om_d) \in (S^1)^d$ satisfies
\[
d ( (\om_1, \om_2, \dots, \om_d),
   \, (\zt_1, \zt_2, \dots, \zt_d) ) < \dt,
\]
then $| \om_j - \zt_j | < \ep$ for $1 \leq j \leq d.$
Then use Lemma~\ref{DensityCond} to choose $n \in \Z^d$ such that
$d ( \sm (n), \, (\zt_1, \zt_2, \dots, \zt_d) ) < \dt.$
Take $v = u_1^{n_1} u_2^{n_2} \cdots u_d^{n_d}$
and use Lemma~\ref{HomToInnW}.
\end{proof}

\begin{lem}\label{MakeSubalgStep1}
Let $\te$ be a nondegenerate skew symmetric real $d \times d$ matrix.
Let $n, \, N \in \N,$ and let $1 \leq k \leq d.$
Then for every $\ep > 0$ there exists
$l = (l_1, l_2, \dots, l_d) \in \Z^d$ such that:
\begin{enumerate}
\item\label{MakeSubalgStep1:1} %
$v = u_1^{l_1} u_2^{l_2} \cdots u_d^{l_d}$ satisfies
$\| v u_j - u_j v \| < \ep$ for $1 \leq j \leq d.$
\item\label{MakeSubalgStep1:2} %
$l_k = 1 \pmod n.$
\item\label{MakeSubalgStep1:3} %
There is $j$ such that $| l_j | > N.$
\end{enumerate}
\end{lem}

\begin{proof}
\Wolog\  $k = 1.$
Set $\af = \Ad (u_1^*).$
There are $\zt_1, \zt_2, \dots, \zt_d \in S^1$ such that
$\af (u_j) = \zt_j u_j$ for $1 \leq j \leq d.$
Let $G = n \Z \oplus \Z^{d - 1},$ which is a
finite index subgroup of $\Z^d.$
According to Corollary~\ref{DensityFinInd}, the subgroup
$\sm (G)$ is dense in $(S^1)^d.$
Let
\[
F = \{ l \in \Z^d
 \colon {\mbox{$| l_j | \leq N + 1$ for $1 \leq j \leq d$}} \}.
\]
Since $F$ is finite, $\sm (G \SM F)$ is also dense in $(S^1)^d.$
Choose $\dt > 0$ small enough that if
$(\om_1, \om_2, \dots, \om_d) \in (S^1)^d$ satisfies
\[
d ( (\om_1, \om_2, \dots, \om_d),
   \, (\zt_1, \zt_2, \dots, \zt_d) ) < \dt,
\]
then $| \om_j - \zt_j | < \ep$ for $1 \leq j \leq d.$
Then use density of $\sm (G \SM F)$ to choose $r \in G \SM F$ such that
$d ( \sm (r), \, (\zt_1, \zt_2, \dots, \zt_d) ) < \dt.$
So with $v_0 = u_1^{r_1} u_2^{r_2} \cdots u_d^{r_d},$
we get $\| v_0 u_j v_0^* - u_1^* u_j u_1 \| < \ep$ for $1 \leq j \leq d.$
Define
\[
l = (r_1 + 1, \, r_2, \, \dots, \, r_d) \in \Z^d
\andeqn
v = u_1^{l_1} u_2^{l_2} \cdots u_d^{l_d} = u_1 v_0 \in A_{\te}.
\]
Clearly $\| v u_j v^* - u_j \| < \ep$ for $1 \leq j \leq d.$
We have $l_1 = 1 \pmod n$ because $r_1 \in n \Z.$
We have $| l_j | > N$ for some $j,$ because
$| r_j | > N + 1$ for some $j.$
\end{proof}

The next lemma is the analog in our context of Lemma~4.6
of~\cite{BKR}.

\begin{lem}\label{MakeSubalgStep2}
Let $\te$ be a nondegenerate skew symmetric real $d \times d$ matrix.
Let $n \in \N,$ let $1 \leq k \leq d,$ and let $\et_0 \in \R \SM \Q.$
Then for every $\ep > 0$ there exist
\[
l = (l_1, l_2, \dots, l_d) \in \Z^d \andeqn
m = (m_1, m_2, \dots, m_d) \in \Z^d
\]
such that:
\begin{enumerate}
\item\label{MakeSubalgStep2:1} %
$v = u_1^{l_1} u_2^{l_2} \cdots u_d^{l_d}$
and $w = u_1^{m_1} u_2^{m_2} \cdots u_d^{m_d}$ satisfy
$\| v u_j - u_j v \| < \ep$ and $\| w u_j - u_j w \| < \ep$
for $1 \leq j \leq d.$
\item\label{MakeSubalgStep2:2} %
There is $\et \in \R \SM \Q$
such that $| \exp ( 2 \pi i \et) - \exp ( 2 \pi i \et_0) | < \ep$
and the unitaries $v$ and $w$ of Part~(\ref{MakeSubalgStep2:1})
satisfy $w v = \exp ( 2 \pi i \et) v w.$
\item\label{MakeSubalgStep2:3} %
$l_k = 1 \pmod n$ and $m_k = 0 \pmod n.$
\end{enumerate}
\end{lem}

\begin{proof}
\Wolog\  $k = 1$ and
$\et_0 \in \left[- \frac{1}{2}, \frac{1}{2} \right].$
Choose $N \in \N$ so large that $2 \pi / N < \ep.$
Use Lemma~\ref{MakeSubalgStep1} with $\te,$ $n,$ and $\ep$ as given,
with $k = 1,$ and with this value of $N,$
obtaining
\[
l \in \Z^d
\andeqn
v = u_1^{l_1} u_2^{l_2} \cdots u_d^{l_d}.
\]
Note in particular that $\| v u_j v^* - u_j \| < \ep$
for $1 \leq j \leq d$ and $l_1 = 1 \pmod n.$
Let $s$ be an index such that $| l_s | > N.$

Let
\[
T = \left\{ \et \in \R \colon
 {\mbox{$\left( u_1^{r_1} u_2^{r_2} \cdots u_d^{r_d} \right) v
        \left( u_1^{r_1} u_2^{r_2} \cdots u_d^{r_d} \right)^*
             = \exp ( 2 \pi i \et) v$ for some $r \in \Z^d$}} \right\}.
\]
Then $T$ is a subgroup of $\R$ which is generated by $d + 1$ elements,
namely $1$ and elements
corresponding to letting $r$ run through the standard basis vectors
of $\Z^d.$
So $T \cap \Q$ is also finitely generated, and is therefore discrete.
Since $\et_0 \not\in \Q,$ we have
$\dist (\et_0, \, T \cap \Q) > 0.$
Set $\ep_0 = \min ( \ep, \, \dist (\et_0, \, T \cap \Q)).$

Set
\[
M = \sum_{j = 1}^d | l_j | \andeqn
\dt = \min \left( \ts{ \frac{1}{2} } \ep_0, \,  M^{-1} \ep_0 \right).
\]
Let $G$ be the finite index subgroup
$G = n \Z \oplus \Z^{d - 1} \S \Z^d.$
Let
\[
\ld = \left(1, \, \dots, \, 1, \, \exp ( 2 \pi i \et_0 / l_s),
 \, 1, \, \dots, \, 1 \right) \in (S^1)^d,
\]
where $\exp ( 2 \pi i \et_0 / l_s)$ is in position $s.$
Use Corollary~\ref{DensityFinInd} and Lemma~\ref{HomToInnW}
to choose $m \in G$ such that
$\sm (m),$ as in Definition~\ref{HomToInnD}, is so close to $\ld$
that $w = u_1^{m_1} u_2^{m_2} \cdots u_d^{m_d}$ satisfies
$\| w u_j w^* - u_j \| < \dt$ for $j \neq s,$ and
$\| w u_s w^* - \exp ( 2 \pi i \et_0 / l_s) u_s \| < \dt.$

Since $\dt \leq \ep,$ it is clear that
$\| w u_j w^* - u_j \| < \ep$ for $j \neq s.$
Also
\[
\| w u_s w^* - u_s \|
 \leq \| w u_s w^* - \exp ( 2 \pi i \et_0 / l_s) u_s \|
   + | \exp ( 2 \pi i \et_0 / l_s) - 1 |.
\]
Using $\dt \leq \frac{1}{2} \ep,$
the first term is less than $\frac{1}{2} \ep.$
The second term satisfies
\[
| \exp ( 2 \pi i \et_0 / l_s) - 1 |
  < 2 \pi \left| \frac{\et_0}{l_s} \right|
  < 2 \pi \left( \frac{1}{2 N} \right)
  \leq {\textstyle{\frac{1}{2}}} \ep.
\]
Therefore $\| w u_j w^* - u_j \| < \ep$ for $j = s$ as well.
This completes the verification of Part~(\ref{MakeSubalgStep2:1})
of the conclusion.
Part~(\ref{MakeSubalgStep2:3}) holds
because $m_1 \in n \Z$ by construction.

It remains to prove Part~(\ref{MakeSubalgStep2:2}).
For each $j$ with $1 \leq j \leq d,$
there is $\zt_j \in S^1$ such that $w u_j w^* = \zt_j u_j.$
Then
\[
w v w^* = \zt_1^{l_1} \zt_2^{l_2} \cdots \zt_d^{l_d} v.
\]
Thus $w v = \exp ( 2 \pi i \et) v w$ for some $\et \in \R.$
By construction we have $| \zt_j - 1 | < M^{-1} \ep_0$ for $j \neq s,$
and $| \zt_s - \exp ( 2 \pi i \et_0 / l_s) | < M^{-1} \ep_0.$
It follows that
\begin{align*}
\left| \rsz{ \zt_1^{l_1} \zt_2^{l_2} \cdots \zt_d^{l_d}
             - \exp ( 2 \pi i \et_0 / l_s)^{l_s} } \right|
  & \leq | l_s | \cdot | \zt_s - \exp ( 2 \pi i \et_0 / l_s) |
      + \sum_{j \neq s} | l_j | \cdot | \zt_j - 1 |  \\
  & < \sum_{j = 1}^d | l_j | M^{-1} \ep_0 \leq \ep_0.
\end{align*}
Therefore
\[
\| w v - \exp ( 2 \pi i \et_0) v w \|
  = \left| \rsz{ \zt_1^{l_1} \zt_2^{l_2} \cdots \zt_d^{l_d}
             - \exp ( 2 \pi i \et_0) } \right|
  < \ep_0,
\]
which is the same as
$| \exp ( 2 \pi i \et) - \exp ( 2 \pi i \et_0) | < \ep_0.$
In particular,
$| \exp ( 2 \pi i \et) - \exp ( 2 \pi i \et_0) | < \ep,$
as desired.
Moreover, $\et \in T$ and there is no $\rh \in T \cap \Q$
such that $| \exp ( 2 \pi i \rh) - \exp ( 2 \pi i \et_0) | < \ep_0,$
whence $\et \not\in \Q.$
\end{proof}

The proofs of the next two results together parallel the proof of
Theorem~1.5 of \cite{BKR}.
The first of them says, roughly,
that higher dimensional noncommutative tori contain approximately
central copies of irrational rotation algebras,
constructed in a special way.
Unfortunately, the rotation parameter varies with the degree
of approximation.

\begin{lem}\label{MakeSubalgStep3}
Let $\te$ be a nondegenerate skew symmetric real $d \times d$ matrix,
let $n \in \N,$ and let $1 \leq k \leq d.$
Then for every $\et_0 \in \R,$
every open set $U \S S^1$ containing $\exp ( 2 \pi i \et_0),$
every finite subset $F \S A_{\te},$
every finite subset $S \S \Gm (U)$
(following Notation~\ref{SectionNot}),
and every $\ep > 0,$
there exist $\et \in \R \SM \Q$ and
\[
l = (l_1, l_2, \dots, l_d) \in \Z^d \andeqn
m = (m_1, m_2, \dots, m_d) \in \Z^d
\]
such that:
\begin{enumerate}
\item\label{MakeSubalgStep3:1} %
$| \et - \et_0 | < \ep$ and $\exp ( 2 \pi i \et) \in U.$
\item\label{MakeSubalgStep3:2} %
$x = u_1^{l_1} u_2^{l_2} \cdots u_d^{l_d}$
and $y = u_1^{m_1} u_2^{m_2} \cdots u_d^{m_d}$ satisfy
$y x = \exp ( 2 \pi i \et) x y.$
\item\label{MakeSubalgStep3:3} %
Following the notation of Theorem~\ref{FieldOfRotAlg},
and with $x$ and $y$ as in Part~(\ref{MakeSubalgStep3:2}),
let $\ph \colon A_{\et} \to A_{\te}$ be the \hm\  such that
$\ph (v_{\et}) = x$ and $\ph (w_{\et}) = y.$
Then $\| [a, \, \ph (b (\et)) ] \| < \ep$
for all $a \in F$ and all $b \in S.$
\item\label{MakeSubalgStep3:4} %
$l_k = 1 \pmod n$ and $m_k = 0 \pmod n.$
\end{enumerate}
\end{lem}

\begin{proof}
Let the notation be as in Theorem~\ref{FieldOfRotAlg}
and Notation~\ref{SectionNot}.

\Wolog\  $\ep < 1.$
Then there is $\ep_0 > 0$ such that whenever $\zt \in S^1$ satisfies
$| \zt - \exp ( 2 \pi i \et_0) | < \ep_0,$
there is a unique $\et \in \R$ such that
$\exp ( 2 \pi i \et) = \zt$ and $| \et - \et_0 | < \ep.$

\Wolog\  $\| a \| \leq 1$ for all $a \in F.$
Replacing $U$ by an open set $V$ with
$\exp ( 2 \pi i \et_0) \in V \S {\overline{V}} \S U,$
we may assume every $b \in S$ is bounded.
Then \wolog\  $\| b (\et) \| \leq 1$ for all $b \in S$ and $\et \in U.$
Write $F = \{ a_1, a_2, \dots, a_s \}$
and $S = \{ b_1, b_2, \dots, b_t \}.$
Choose polynomials $g_1, g_2, \dots, g_t$
in four noncommuting variables such that
\[
\| g_r (v_{\et_0}, \, v_{\et_0}^*, \, w_{\et_0}, \, w_{\et_0}^*)
        - b_r (\et_0) \|
   < \ts{ \frac{1}{7} \ep }
\]
for $1 \leq r \leq t.$
Because the rotation algebras form a \ct\  field over $S^1$
(Theorem~\ref{FieldOfRotAlg}), there is $\dt > 0$ such that
whenever $| \et - \et_0 | < \dt$ we have $\exp ( 2 \pi i \et_0) \in U,$
and
\[
\| g_r (v_{\et}, \, v_{\et}^*, \, w_{\et}, \, w_{\et}^*) - b_r (\et) \|
   < \ts{ \frac{2}{7} \ep }
\]
for $1 \leq r \leq t.$

Choose polynomials $f_1, f_2, \dots, f_t$
in $2 d$ noncommuting variables such that
\[
\| f_r (u_1, \, u_1^*, \, \dots, \, u_d, \, u_d^*) - a_r \|
   < \frac{\ep}{7 (1 + \ep)}
\]
for $1 \leq r \leq s.$
Choose (see Proposition~4.3 of~\cite{BKR}) $\dt_0 > 0$ such that
whenever $D$ is a \ca\  and
\[
c_1, c_2, \dots, c_{2 d}, d_1, d_2, d_3, d_4 \in D
\]
are elements of norm~$1$
which satisfy $\| [c_r, d_j] \| < \dt_0$ for all $j$ and $r,$ then
\[
\| [f_r (c_1, c_2, \dots, c_{2 d}), \,
               g_j (d_1, d_2, d_3, d_4)] \|
  < \ts{ \frac{1}{7} \ep }
\]
for $1 \leq r \leq s$ and $1 \leq j \leq t.$

Apply Lemma~\ref{MakeSubalgStep2} with $\te,$ $n,$ $\et_0,$ and $k$
as given,
and with $\min (\ep_0, \dt, \dt_0)$ in place of $\ep.$
We obtain $\et \in \R \SM \Q$ and
\[
l = (l_1, l_2, \dots, l_d) \in \Z^d \andeqn
m = (m_1, m_2, \dots, m_d) \in \Z^d.
\]
Set
\[
x = u_1^{l_1} u_2^{l_2} \cdots u_d^{l_d} \andeqn
y = u_1^{m_1} u_2^{m_2} \cdots u_d^{m_d}.
\]
By the choice of $\ep_0,$ we may assume that $| \et - \et_0 | < \ep,$
and by the choice of $\dt$ we have $\exp ( 2 \pi i \et_0) \in U.$
This is Part~(\ref{MakeSubalgStep3:1}) of the conclusion.
Parts~(\ref{MakeSubalgStep3:2}) and~(\ref{MakeSubalgStep3:4})
are immediate.

It remains to prove Part~(\ref{MakeSubalgStep3:3}).
Part~(\ref{MakeSubalgStep2:1})
of the conclusion of Lemma~\ref{MakeSubalgStep2} and
the choice of $\dt_0$ ensure that
\[
\| [f_r (u_1, \, u_1^*, \, \dots, \, u_d, \, u_d^*), \,
               g_j (x, x^*, y, y^*)] \|
  < \ts{ \frac{1}{7} \ep }
\]
for $1 \leq r \leq s$ and $1 \leq j \leq t.$
{}From the choice of $\dt,$ we get
\[
\| g_j (x, x^*, y, y^*) \|
   < \| \ph (b_j (\et)) \| + \ts{ \frac{2}{7} \ep }
   < 1 + \ep
\]
for $1 \leq j \leq t.$
Using the choice of the polynomials $f_r,$ we therefore get
\begin{align*}
\| [a_r, \, \ph (b_j (\et)) ] \|
& \leq 2 \| a_r \| \cdot \| \ph (b_j (\et))
            - g_j (x, x^*, y, y^*) \|  \\
& \hspace*{3em} \mbox{}
     + 2 \| a_r - f_r (u_1, \, u_1^*, \, \dots, \, u_d, \, u_d^*) \|
               \cdot \| g_j (x, x^*, y, y^*) \|  \\
& \hspace*{3em} \mbox{}
     + \| [f_r (u_1, \, u_1^*, \, \dots, \, u_d, \, u_d^*), \,
               g_j (x, x^*, y, y^*)] \|  \\
& < 2 \left( \frac{2 \ep}{7} \right)
     + 2 (1 + \ep) \left( \frac{\ep}{7 (1 + \ep)} \right)
     + \frac{\ep}{7} = \ep
\end{align*}
for $1 \leq r \leq s$ and $1 \leq j \leq t,$
as desired.
\end{proof}

\begin{prp}\label{IrratAPR}
Let $\te$ be a nondegenerate skew symmetric real $d \times d$ matrix.
Let $n \in \N,$ let $\om = \exp (2 \pi i / n),$
let $1 \leq k \leq d,$ and, following Notation~\ref{StdNtn}, let
$\af \colon A_{\te} \to A_{\te}$ the unique automorphism
satisfying $\af (u_k) = \om u_k$ and $\af (u_r) = u_r$
for $r \neq k.$
Then the action of $\Zqn$
generated by $\af$ has the tracial Rokhlin property.
\end{prp}

\begin{proof}
Let $\ta$ be the unique tracial state on $A_{\te}$
(Theorem~\ref{Simplicity}).
We will show that for every $\ep > 0$ and
every finite subset $F \S A_{\te},$
there are
\mops\  $e_0, e_1, \dots, e_{n - 1} \in A_{\te}$ such that:
\begin{enumerate}
\item\label{IrratAPR:1} %
$\| \af (e_j) - e_{j + 1} \| < \ep$ for $0 \leq j \leq n - 1.$
\item\label{IrratAPR:2} %
$\| e_j a - a e_j \| < \ep$
for $0 \leq j \leq n - 1$ and $a \in F.$
\item\label{IrratAPR:3} %
$1 - n \ta (e_0) < \ep.$
\end{enumerate}

We first argue that this is enough to deduce the \tRp.
We must prove Condition~(\ref{TRPCond:3}) in Proposition~\ref{TRPCond}.
We first recall
(see Theorems~6.1 and~7.1 of~\cite{Rf1},
or Theorems~1.4(d) and~1.5 of~\cite{BKR})
that if $p, \, q \in A_{\te}$ are \pj s with $\ta (p) < \ta (q),$
then $p \precsim q.$
Also, $A_{\te}$ has Property~(SP) by Theorem~1.4(b) of~\cite{BKR}.
If now a nonzero positive element $x \in A_{\te}$ is given,
then we may use Property~(SP)
to find a nonzero \pj\  $p \in {\overline{x A x}}.$
Require $\ep \leq \min (\ta (p), 1).$
Then $\af (e_j) \sim e_{j + 1}.$
Let $e = \sum_{j = 0}^{n - 1} e_j.$
This gives $\ta (1 - e) = 1 - n \ta (e_0) < \ep,$
whence $\ta (e_0) > \frac{1}{n} (1 - \ep).$
Now $\ta (1 - e) < \ep,$ which implies
$\ta (1 - e) < \ta (p),$
so that Condition~(\ref{TRPCond:3}) of Proposition~\ref{TRPCond}
follows from the comparison result above.

Now we prove Conditions (\ref{IrratAPR:1}), (\ref{IrratAPR:2}),
and~(\ref{IrratAPR:3}) at the beginning of the proof.
Let the notation be as in Theorem~\ref{FieldOfRotAlg}
and Notation~\ref{SectionNot}.
Let $\ep > 0.$
Choose and fix $\et_0 \in \R \SM \Q.$
Choose $\ep_1 > 0$ such that whenever $a_0, a_1, \dots, a_{n - 1}$
are elements of a unital \ca\  $D$ with
\[
\| a_j a_r - \dt_{j, r} a_j \| < \ep_1 \andeqn
\| a_j^* - a_j \| < \ep_1
\]
for $0 \leq j, \, r \leq n - 1,$ then there are \mops\  %
\[
q_0, q_1, \dots, q_{n - 1} \in D
\]
such that
$\| q_j - a_j \| < \frac{1}{3} n^{-1} \ep$ for $0 \leq j \leq n - 1.$
(For example, apply Definition~2.2 and Lemma~2.3 of~\cite{BKR}
with the \fd\  \ca\  $B$ taken to be $\C^{n + 1},$ using in
addition the element $a_n = 1 - \sum_{j = 0}^{n - 1} a_j.$)
Let $p_0, p_1, \dots, p_{n - 1} \in A_{\et_0}$ be the \pj s
$e_0, e_1, \dots, e_{n - 1}$ of Lemma~\ref{2D}
for $\et_0$ in place of $\et$ and $\frac{1}{3} \ep$ in place of $\ep.$
Because the rotation algebras form a \ct\  field over $S^1$
with section algebra $A$ (Theorem~\ref{FieldOfRotAlg}),
we may choose $c_0, c_1, \dots, c_{n - 1} \in A$
such that $\ev_{\et_0} (c_j) = p_j$ for $0 \leq j \leq n - 1,$
and we can furthermore find $\dt_0 > 0$ such that
$| \exp ( 2 \pi i \et) - \exp ( 2 \pi i \et_0) | < \dt_0$ implies
\[
\| \ev_{\et} (c_j) \ev_{\et} (c_r)
   - \dt_{j, r} \ev_{\et} (c_j) \| < \ep_1
\andeqn
\| \ev_{\et} (c_j)^* - \ev_{\et} (c_j) \| < \ep_1
\]
for $0 \leq j, \, r \leq n - 1.$
Let $V \S S^1$ be an open set such that $\exp ( 2 \pi i \et_0) \in V$
and such that $\zt \in {\overline{V}}$ implies
$| \zt - \exp ( 2 \pi i \et_0) | < \dt_0.$
Letting $c_j |_{\overline{V}}$ denote the restriction of $c_j,$
regarded as a section, to ${\overline{V}},$
we get
\[
\left\| \ts{ \left( c_j |_{\overline{V}} \right) }
               \ts{ \left( c_r |_{\overline{V}} \right) }
        - \dt_{j, r} c_j |_{\overline{V}} \right\|
     < \ep_1 \andeqn
\left\| \ts{ \left( c_j |_{\overline{V}} \right) }^*
                - c_j |_{\overline{V}} \right\| < \ep_1
\]
for $0 \leq j, \, r \leq n - 1,$
so that there are \mops\  %
\[
q_0, q_1, \dots, q_{n - 1} \in \Gm ( {\overline{V}} )
\]
such that
$\left\| q_j - c_j |_{\overline{V}} \right\| < \frac{1}{3} n^{-1} \ep$
for $0 \leq j \leq n - 1.$
Since the restriction map $A = \Gm (S^1) \to \Gm ( {\overline{V}} )$
is surjective,
there exist $b_0, b_1, \dots, b_{n - 1} \in A$
such that $b_j |_{\overline{V}} = q_j$ for $0 \leq j \leq n - 1.$

Let the generators of $A$ be as in Theorem~\ref{FieldOfRotAlg},
and let $\bt \in \Aut (A)$ be the unique automorphism such that
\[
\bt (v) = \om v, \,\,\,\,\,\, \bt (w) = w, \andeqn \bt (z) = z.
\]
Let $\bt_{\et} \in \Aut (A_{\et})$ be defined by
$\bt_{\et} ( v_{\et}) = \om v_{\et}$
and $\bt_{\et} ( w_{\et}) = w_{\et}.$
Then $\ev_{\et} \circ \bt = \bt_{\et} \circ \ev_{\et}.$
Since $\bt$ sends continuous sections to continuous sections,
there is an open set $U_0 \S V$ such that $\et_0 \in U_0$
and if $\et \in U_0$
then for $0 \leq j \leq n - 1$ and with $b_n = b_0,$
\[
\| \bt_{\et} ( b_j (\et) ) - b_{j + 1} (\et) \|
\andeqn
\| \bt_{\et_0} ( b_j (\et_0) ) - b_{j + 1} (\et_0) \|
\]
differ by less than $\frac{1}{3} \ep.$
For such $\et$ we have $b_j (\et) = q_j (\et),$ so,
using $c_j (\et_0) = p_j$ and $\bt_{\et_0} ( p_j) = p_{j + 1}$
at the second last step,
\begin{align*}
\| \bt_{\et} ( q_j (\et) ) - q_{j + 1} (\et) \|
& < \| \bt_{\et_0} ( q_j (\et_0) ) - q_{j + 1} (\et_0) \|
        + \ts{ \frac{1}{3} } \ep     \\
& \hspace*{-2em}
  < \left\| q_j - c_j |_{\overline{V}} \right\|
        + \left\| q_{j + 1} - c_{j + 1} |_{\overline{V}} \right\|
        + \| \bt_{\et_0} ( c_j (\et_0) ) - c_{j + 1} (\et_0) \|
        + \ts{ \frac{1}{3} } \ep     \\
& \hspace*{-2em}
  < \ts{ \frac{1}{3} } n^{-1} \ep + \ts{ \frac{1}{3} } n^{-1} \ep
        + \ts{ \frac{1}{3} } \ep \leq \ep
\end{align*}
for $0 \leq j \leq n - 1.$

Using Lemma~\ref{ContOfTrace},
choose an open set $U \S U_0$ such that $\et_0 \in U$
and if $\et \in U$ then for $0 \leq j \leq n - 1$ we have
$| \ta_{\et} (q_j (\et)) - \ta_{\et_0} (q_j (\et_0)) |
                < \ts{ \frac{1}{3} } n^{-1} \ep.$

Apply Lemma~\ref{MakeSubalgStep3} with $\te,$ $n,$ $k,$ $\et_0,$
$U,$ and $F$ as given,
with $\min ( \ep, \dt)$ in place of $\ep,$
and with $S = \{ q_0, \, q_1, \, \dots, \, q_{n - 1} \}.$
We obtain $\et \in (\R \SM \Q) \cap U$ and
\[
l = (l_1, l_2, \dots, l_d) \in \Z^d \andeqn
m = (m_1, m_2, \dots, m_d) \in \Z^d.
\]
Set
\[
x = u_1^{l_1} u_2^{l_2} \cdots u_d^{l_d} \andeqn
y = u_1^{m_1} u_2^{m_2} \cdots u_d^{m_d},
\]
so that $y x = \exp ( 2 \pi i \et) x y.$
Let $\ph \colon A_{\et} \to A_{\te}$ be the \hm\  such that
$\ph (v_{\et}) = x$ and $\ph (w_{\et}) = y,$
and set $e_j = \ph ( q_j (\et))$ for $0 \leq j \leq n - 1.$
We verify Conditions~(\ref{IrratAPR:1}), (\ref{IrratAPR:2}),
and~(\ref{IrratAPR:3}) at the beginning of the
proof for this choice of $e_0, e_1, \dots, e_{n - 1}.$

We do Condition~(\ref{IrratAPR:1}).
Because $l_k = 1 \pmod n$ and $m_k = 0 \pmod n,$ we have
$\af (x) = \om x$ and $\af (y) = y.$
It follows that $\af \circ \ph = \ph \circ \bt_{\et}.$
Therefore
\[
\| \af (e_j) - e_{j + 1} \|
   \leq \| \bt_{\et} ( q_j (\et) ) - q_{j + 1} (\et) \| < \ep
\]
for $0 \leq j \leq n - 1,$ as desired.

Condition~(\ref{IrratAPR:2}) is immediate from
Part~(\ref{MakeSubalgStep3:3}) of Lemma~\ref{MakeSubalgStep3}.

Finally, we check Condition~(\ref{IrratAPR:3}).
By uniqueness of the tracial states, we have
$\ta \circ \ph = \ta_{\et}.$
Therefore, using the choice of $U$ at the second step and
$\| q_j (\et_0) - p_j \| < \ts{ \frac{1}{3} } n^{-1} \ep$
at the third step, we get
\[
\ta (e_j)
   = \ta_{\et} (q_j (\et))
   > \ta_{\et_0} (q_j (\et_0)) - \ts{ \frac{1}{3} } n^{-1} \ep
   > \ta_{\et_0} (p_j) - \ts{ \frac{2}{3} } n^{-1} \ep.
\]
Therefore
\[
1 - n \ta (e_0)
  < 1 - n \ta (p_0) + \ts{ \frac{2}{3} } \ep
  < \ts{ \frac{1}{3} } \ep + \ts{ \frac{2}{3} } \ep = \ep.
\]
This completes the proof of~(\ref{IrratAPR:3}).
\end{proof}

We next need to identify the fixed point algebra of the action
in Proposition~\ref{IrratAPR} with a suitable
higher dimensional noncommutative torus.
The following lemma, suggested by Hanfeng Li,
is a substantial generalization of our original statement.

\begin{lem}\label{FPOfMult}
Let $\te$ be a skew symmetric real $d \times d$ matrix.
Let
\[
M = ( m_{j, k} )_{1 \leq j, k \leq d}
           \in {\mathrm{GL}}_d (\R) \cap M_n (\Z),
\]
and set ${\widetilde{\te}} = M^{\mathrm{t}} \te M.$
Let $u_1, u_2, \ldots, u_d$ be the standard generators
of $A_{\te}$ (as in Notation~\ref{StdNtn}),
and let
${\widetilde{u}}_1, {\widetilde{u}}_2, \ldots, {\widetilde{u}}_d$
be the standard generators of $A_{\widetilde{\te}}.$
For $1 \leq k \leq d,$ define
\[
v_k = u_1^{m_{1, k}} \cdot u_2^{m_{2, k}} \cdots u_d^{m_{d, k}}
  \in A_{\te}.
\]
Then ${\widetilde{u}}_k \mapsto v_k$ extends to an isomorphism
$A_{\widetilde{\te}} \to C^* (v_1, v_2, \ldots, v_d).$
\end{lem}

\begin{proof}
It is easy to check that there exists
a homomorphism $\ph \colon A_{\widetilde{\te}} \to A_{\te}$
such that $\ph \big( {\widetilde{u}}_k \big) = v_k$
for $1 \leq k \leq d,$
and clearly $\ph ( A_{\widetilde{\te}} ) = C^* (v_1, v_2, \ldots, v_d).$
We need only check that $\ph$ is injective.

Define a group action $\gm \colon (S^1)^d \to \Aut (A_{\te})$
by $\gm_{\zt_1, \zt_2, \ldots, \zt_d} (u_k) = \zt_k u_k$
for $1 \leq k \leq d$ and $\zt_1, \zt_2, \ldots, \zt_d \in S^1.$
Similarly define
${\widetilde{\gm}} \colon (S^1)^d \to \Aut (A_{\widetilde{\te}}).$
It is well known (and is easily checked by averaging monomials
in the $u_k$ and ${\widetilde{u}}_k$ over the group)
that the fixed point algebras $A_{\te}^{\gm}$
and $A_{\widetilde{\te}}^{\widetilde{\gm}}$ are both $\C \cdot 1.$
In particular,
the restriction of $\ph$ to $A_{\widetilde{\te}}^{\widetilde{\gm}}$
is injective.

Since $M \in M_d (\Z),$ we have $M^{\mathrm{t}} \Z^d \subset \Z^d,$
so that $M^{\mathrm{t}}$ descends to a \hm\  %
$f \colon (S^1)^d \to (S^1)^d.$
This \hm\  is surjective because $M$ is invertible over $\R.$
Define $\bt \colon (S^1)^d \to \Aut (A_{\widetilde{\te}})$
by
$\bt_{\zt_1, \zt_2, \ldots, \zt_d}
       = {\widetilde{\gm}}_{f (\zt_1, \zt_2, \ldots, \zt_d)}.$
Then one checks,
by examining the generators ${\widetilde{u}}_k,$
that
$\ph \circ \bt_{\zt_1, \zt_2, \ldots, \zt_d}
   = \gm_{\zt_1, \zt_2, \ldots, \zt_d} \circ \ph$
for $\zt_1, \zt_2, \ldots, \zt_d \in S^1.$
Injectivity of $\ph |_{A_{\widetilde{\te}}^{\widetilde{\gm}} }$
now implies injectivity of $\ph,$
since if $a$ is a nonzero positive element of ${\mathrm{Ker}} (\ph),$
then averaging over $(S^1)^d$ gives a
nonzero positive element of
${\mathrm{Ker}} (\ph) \cap A_{\widetilde{\te}}^{\widetilde{\gm}}.$
\end{proof}

\begin{cor}\label{IrrFPIsTAF}
Let $\te$ be a nondegenerate skew symmetric real $d \times d$ matrix.
Let $n \in \N,$ let $1 \leq l \leq d,$ and let
\[
B = \diag (1, \dots, 1, n, 1, \dots, 1) \in {\mathrm{GL}}_d (\Q),
\]
where $n$ is in the $l$-th position.
Then $A_{B^{\mathrm{t}} \te B}$ has tracial rank zero
\ifo\  $A_{\te}$ has tracial rank zero.
\end{cor}

\begin{proof}
With $u_k$ as in Notation~\ref{StdNtn}, set
\[
D = C^* (u_1, \, \dots, \, u_{l - 1}, \, u_l^n,
           \, u_{l + 1}, \, \dots, \, u_d)
  \S A_{\te}.
\]
Calculating $B^{\mathrm{t}} \te B,$
we find that $D \cong A_{B^{\mathrm{t}} \te B}$ by Lemma~\ref{FPOfMult}.

Let
$\af \colon A_{\te} \to A_{\te}$ the unique automorphism
satisfying $\af (u_l) = \exp (2 \pi i / n) u_l$
and $\af (u_k) = u_k$ for $k \neq l.$
We claim that $A_{\te}^{\af} = D.$
That $D \S A_{\te}^{\af}$ is clear.
For the reverse inclusion, define $E \colon A_{\te} \to A_{\te}^{\af}$
by $E (a) = \frac{1}{n} \sum_{j = 0}^{n - 1} \af^j (a).$
Then $E$ is a surjective continuous linear map, so it suffices to
show that $E \big( u_1^{m_1} u_2^{m_2} \cdots u_d^{m_d} \big) \in D$
for all $m = (m_1, m_2, \dots, m_d) \in \Z^d.$
If $m_l$ is divisible by $n$ then
$u_1^{m_1} u_2^{m_2} \cdots u_d^{m_d}$
is a fixed point of $E$ and is in $D,$
and otherwise
$E \big( u_1^{m_1} u_2^{m_2} \cdots u_d^{m_d} \big) = 0 \in D.$
This proves the claim.

Using
Proposition~\ref{IrratAPR} and Corollary~\ref{GaugeAppInn},
the result now follows from Theorem~4.8 of~\cite{PhtRp1a}.
\end{proof}

\section{Direct limit decomposition for simple noncommutative
            tori}\label{Sec:NCTInd}

\indent
In this section, we use the results of the previous two sections
to prove that every simple higher dimensional noncommutative torus
is an AT~algebra.

The following result is essentially Corollary~6.6 of~\cite{Ks4}.

\begin{prp}\label{SingleStep}
Let $\af$ be a nondegenerate skew symmetric real bicharacter on $\Z^n.$
Suppose that $A_{\af |_{\Z^{n - 1} \times \{ 0 \} }}$ is a
simple AT~algebra with real rank zero.
Then $A_{\af}$ is a simple AT~algebra with real rank zero.
\end{prp}

\begin{proof}
Let $\bt = \af |_{\Z^{n - 1} \times \{ 0 \} }.$
We note that
$K_0 (A_{\bt}) \cong K_1 (A_{\bt}) \cong \Z^{2^{n - 1}}$
by Lemma~\ref{ItCrPrd} and by repeated application of the
Pimsner-Voiculescu exact sequence~\cite{PV}.
In particular, both groups are finitely generated.
Further write $A_{\af} = C^* ( \Z, A_{\bt}, \ph)$ as in
Lemma~\ref{ItCrPrd}, with $\ph$ homotopic to the identity.
Thus, in the notation of~\cite{Ks4}
(see the introduction to~\cite{Ks4}),
$\ph \in {\mathrm{HInn}} (A_{\bt}).$
So the proof of Corollary~6.5 of~\cite{Ks4} shows that the
hypotheses of Theorem~6.4 of~\cite{Ks4} hold.
We know from Lemma~\ref{Simplicity} that
$A_{\af} = C^* ( \Z, A_{\bt}, \ph)$ has a unique tracial state.
Therefore Theorem~6.4 of~\cite{Ks4} implies that
$A_{\af} = C^* ( \Z, A_{\bt}, \ph)$
is a simple AT~algebra with real rank zero.
\end{proof}

It is worth pointing out another version of the proof.
The automorphism $\ph$ has the Rokhlin property
by Theorem~6.1 of~\cite{Ks2} and the preceding remark in~\cite{Ks2}.
The algebra $A_{\bt}$ has tracial rank zero,
since it is a simple unital AT~algebra with real rank zero.
Therefore $\ph$ has the \tRp\  for actions of $\Z,$
by Theorem~1.12 of~\cite{OP1}.
It is approximately inner by Corollary~\ref{GaugeAppInn}.
So Theorem~3.9 of~\cite{LO} applies,
showing that $A_{\af}$ has tracial rank zero.
Now one can conclude that $A_{\af}$ is an AT~algebra
by the same proof as for Theorem~\ref{NCTIsAT}.
However, tracial rank zero is all that is really needed here.

\begin{lem}\label{GenOfGLQ}
The group ${\mathrm{GL}}_d (\Q)$ is generated as a group by
${\mathrm{GL}}_d (\Z)$ and all matrices of the form
$\diag (1, \dots, 1, n, 1, \dots, 1),$ where $n \in \N$ is nonzero
and is in an arbitrary position.
\end{lem}

\begin{proof}
Let $G$ be the subgroup of ${\mathrm{GL}}_d (\Q)$ generated
by ${\mathrm{GL}}_d (\Z)$ and the matrices
$\diag (1, \dots, 1, n, 1, \dots, 1).$
It suffices to show that $G$ contains all of the following three
kinds of elementary matrices:
\[
E_j^{(1)} (r) = \diag (1, \dots, 1, r, 1, \dots, 1),
\]
where $r \in \Q \SM \{ 0 \}$ and is the $j$-th diagonal entry in the
matrix;
the transposition matrix $E_{j, k}^{(2)},$ for $1 \leq j < k \leq d,$
which acts on the standard basis vectors by
\[
E_{j, k}^{(2)} (\dt_l) = \left\{ \begin{array}{ll}
     \dt_l   & \hspace{3em}  l \neq j, \, k \\
     \dt_k   & \hspace{3em}  l = j  \\
     \dt_j   & \hspace{3em}  l = k;
    \end{array} \right.
\]
and the matrix $E_{j, k}^{(3)} (r)$
for $1 \leq j, \, k \leq n$ with $j \neq k$ and $r \in \Q,$
given by
\[
E_{j, k}^{(3)} (r) (\dt_l) = \left\{ \begin{array}{ll}
     \dt_l             & \hspace{3em}  l \neq k \\
     \dt_k + r \dt_j   & \hspace{3em}  l = k.
    \end{array} \right.
\]

If $r = (-1)^m p / q$ with $m = 0$ or $m = 1$ and with $p$ and $q$
positive integers, then
\[
E_j^{(1)} (r)
 = E_j^{(1)} ((-1)^m) E_j^{(1)} (p) \big[ E_j^{(1)} (q) \big]^{-1},
\]
where the first factor is in ${\mathrm{GL}}_d (\Z)$ and
$E_j^{(1)} (p)$ and $E_j^{(1)} (q)$ are among the other generators
of $G.$
The matrix $E_{j, k}^{(2)}$ is already in ${\mathrm{GL}}_d (\Z).$
For $E_{j, k}^{(3)} (r),$ we may conjugate by a permutation matrix,
which is in ${\mathrm{GL}}_d (\Z),$ and split off as a direct summand
a $(d - 2) \times (d - 2)$ identity matrix, and thus reduce to the
case $d = 2,$ $j = 1,$ and $k = 2.$
Write $r = p / q$ with $p \in \Z$ and $q \in \N.$
Then the factorization
\[
E_{1, 2}^{(3)} (r)
 = \left( \begin{array}{cc} 1 & p/q \\ 0 & 1 \end{array} \right)
 = \left( \begin{array}{cc} q & 0 \\ 0 & 1 \end{array} \right)^{-1}
   \left( \begin{array}{cc} 1 & p \\ 0 & 1 \end{array} \right)
   \left( \begin{array}{cc} q & 0 \\ 0 & 1 \end{array} \right)
\]
shows that $E_{1, 2}^{(3)} (r) \in G.$
\end{proof}

\begin{cor}\label{TAFConjByRat}
Let $\te$ be a nondegenerate skew symmetric real $d \times d$ matrix.
Let $B \in {\mathrm{GL}}_d (\Q).$
Then $A_{B^{\mathrm{t}} \te B}$ has tracial rank zero
\ifo\  $A_{\te}$ has tracial rank zero.
\end{cor}

\begin{proof}
Combine Lemma~\ref{GenOfGLQ}, Remark~\ref{CoordFree}, and
Corollary~\ref{IrrFPIsTAF}.
\end{proof}

\begin{lem}\label{StrOfNondeg}
Let $\te$ be a nondegenerate skew symmetric real $d \times d$ matrix,
with $d > 2.$
Suppose that there is no subgroup $H$ of $\Z^d$ of rank $d - 1$
such that $\te |_H$
(in the sense of Remark~\ref{Restrict}) is nondegenerate.
Let $r < d - 1$ be the maximal rank of a proper subgroup $H$ of $\Z^d$
such that $\te |_H$ is nondegenerate.
Then there exists $B \in {\mathrm{GL}}_d (\Q)$ such that
$B^{\mathrm{t}} \te B$ has the block form
\[
B^{\mathrm{t}} \te B
   = \left( \begin{array}{cc} \rh_{1, 1} & \rh_{1, 2} \\
           (\rh_{1, 2})^{\mathrm{t}} & \rh_{2, 2} \end{array} \right),
\]
with $\rh_{1, 1}$ and $\rh_{2, 2}$
nondegenerate skew symmetric real $r \times r$ and
$(d - r) \times (d - r)$ matrices, and where all the entries of
$\rh_{1, 2}$ are in $\Z.$
\end{lem}

\begin{proof}
Let $H \S \Z^d$ be a subgroup of rank $r$ such that
$\te |_H$ is nondegenerate.
Since $\te$ is not rational, clearly $r \geq 2.$
Let $(v_1, v_2, \dots, v_r)$ be a basis for $H$ over $\Z.$
Choose $v_{r + 1}, \, \dots, \, v_d \in \Z^d$ such that
$(v_1, v_2, \dots, v_d)$ is a basis for $\Q^d$ over $\Q.$
For $r + 1 \leq k \leq d,$ by hypothesis
$\te |_{H + \Z v_k}$ is degenerate.
By Lemma~\ref{NondegCond},
there exists $x_k \in \spn_{\Q} (H \cup \{ v_k \} ) \SM \{ 0 \}$
such that $\langle x_k, \, \te y \rangle \in \Q$
for all $y \in \spn_{\Q} (H \cup \{ v_k \} ).$
Since $\te |_H$ is nondegenerate, we have $x_k \not\in \spn_{\Q} (H).$
Therefore $v_k \in \spn_{\Q} (H \cup \{ x_k \} ).$
It follows that
\[
v_{r + 1}, \, \dots, \, v_d \in \spn_{\Q}
 (v_1, \, v_2, \, \dots, \, v_r, \, x_{r + 1}, \, \dots, \, x_d),
\]
so that
$(v_1, \, v_2, \, \dots, \, v_r, \, x_{r + 1}, \, \dots, \, x_d)$
is a basis for $\Q^d.$
By construction, we have
$\langle x_k, \te v_l \rangle \in \Q$ for $1 \leq l \leq r$ and
$r + 1 \leq k \leq d.$
Choose $N \in \Z \SM \{ 0 \}$ such that
$N \langle x_k, \te v_l \rangle \in \Z$ for $1 \leq l \leq r$ and
$r + 1 \leq k \leq d.$

Let $B \in {\mathrm{GL}}_d (\Q)$ be the matrix whose
action on the standard basis vectors is
\[
B \dt_k = \left\{ \begin{array}{ll}
     v_k   & \hspace{3em}  1 \leq k \leq r  \\
    N x_k  & \hspace{3em}  r + 1 \leq k \leq d.
    \end{array} \right.
\]
Then for $1 \leq l \leq r$ and $r + 1 \leq k \leq d,$ we have
\[
\langle \dt_k, B^{\mathrm{t}} \te B \dt_l \rangle
   = N \langle x_k, \te v_l \rangle \in \Z.
\]
Since $B^{\mathrm{t}} \te B$ is skew symmetric, this shows that it
has a block decomposition of the required form.
It is immediate to check that the two diagonal blocks must be
nondegenerate, since otherwise $B^{\mathrm{t}} \te B$ would be
degenerate, contradicting Lemma~\ref{ConjByRat}.
\end{proof}

\begin{thm}\label{NCTIsTAF}
Let $\te$ be a nondegenerate skew symmetric real $d \times d$ matrix,
with $d \geq 2.$
Then $A_{\te}$ has tracial rank zero.
\end{thm}

\begin{proof}
We prove this by induction on $d.$
For $d = 2,$ the Elliott-Evans Theorem~\cite{EE}
shows that $A_{\te}$ is a simple AT~algebra with real rank zero,
and tracial rank zero then follows from Proposition~2.6 of~\cite{LnTAF}
(with ${\mathcal{C}}$ as defined in~2.5 of~\cite{LnTAF}).
Suppose $d$ is given, and the theorem is known for all
skew symmetric real $k \times k$ matrices with $k < d.$
Let $\te$ be a nondegenerate skew symmetric real $d \times d$ matrix.
There are two cases.

First, suppose that there is a subgroup $H_0$ of $\Z^d$ of rank $d - 1$
such that $\te |_{H_0}$
(in the sense of Remark~\ref{Restrict}) is nondegenerate.
Set
\[
H = \{ x \in \Z^d \colon
 {\mbox{There is $n \in \Z$ such that $n x \in H_0$}} \}.
\]
Then $H$ is also a subgroup of $\Z^d$ of rank $d - 1,$ and
$\te |_{H_0}$ is also nondegenerate.
Moreover, $\Z^d / H$ is torsion free and therefore isomorphic to $\Z,$
from which it follows that the quotient map splits.
Thus there is an isomorphism $\Z^d \to \Z^d$ which sends
$H$ isomorphically onto $\Z^{d - 1} \oplus \{ 0 \} \S \Z^d.$
Accordingly, we may assume that $H = \Z^{d - 1} \oplus \{ 0 \}.$
By the induction hypothesis,
$A_{\te |_H}$ is a simple AT~algebra with real rank zero.
So Proposition~\ref{SingleStep} implies that
$A_{\te}$ is a simple AT~algebra with real rank zero.

Now assume there is no such subgroup $H_0$ of rank $d - 1.$
Let $B$ be as in Lemma~\ref{StrOfNondeg}, with
\[
B^{\mathrm{t}} \te B
   = \left( \begin{array}{cc} \rh_{1, 1} & \rh_{1, 2} \\
             (\rh_{1, 2})^{\mathrm{t}} & \rh_{2, 2} \end{array} \right),
\]
and where in particular all the entries of
$\rh_{1, 2}$ are in $\Z.$
Then $A_{B^{\mathrm{t}} \te B} \cong A_{\rh}$ for
\[
\rh = \left( \begin{array}{cc} \rh_{1, 1} & 0 \\
                                   0   & \rh_{2, 2} \end{array} \right).
\]
Using the definitions of $A_{\rh},$ $A_{ \rh_{1, 1} },$ and
$ A_{ \rh_{2, 2} }$ as universal algebras on generators and relations,
one easily checks that
$A_{\rh} \cong A_{ \rh_{1, 1} } \otimes A_{ \rh_{2, 2} }.$
By the induction hypothesis,
both $A_{ \rh_{1, 1} }$ and $A_{ \rh_{2, 2} }$
are simple AT~algebras with real rank zero.
Therefore $A_{ \rh_{1, 1} } \otimes A_{ \rh_{2, 2} }$ is a simple
direct limit, with no dimension growth, of homogeneous \ca s.
Since it has a unique tracial state, Theorems~1 and~2 of~\cite{BDR}
imply that $A_{\rh}$ has stable rank one and real rank zero.
This algebra has weakly unperforated
K-theory by Theorem~6.1 of~\cite{Rf1}.
(Actually, this is true for any direct limit of the type at hand.)
It now follows from Theorem~4.6 of~\cite{Ln5} that
$A_{ \rh_{1, 1} } \otimes A_{ \rh_{2, 2} }$ has tracial rank zero.
So Corollary~\ref{TAFConjByRat}
shows that $A_{\te}$ has tracial rank zero.
\end{proof}

We note that one could use the earlier Theorem~3.11 of~\cite{EG}
to show that $A_{ \rh_{1, 1} } \otimes A_{ \rh_{2, 2} }$ is an
AT~algebra with real rank zero, from which it follows that
this algebra has tracial rank zero.
The use of H.~Lin's classification theorem,
Theorem~5.2 of~\cite{Ln15},
remains essential,
because we can only prove that 
crossed products and fixed point algebras of actions by
finite cyclic groups with the tracial Rokhlin property
preserve tracial rank zero,
not that they the property of being an AT~algebra
or even an AH~algebra.

\begin{rmk}\label{NoHg}
Since the paper~\cite{Hg} remains unpublished, it is worth pointing
out that the proof of Theorem~\ref{NCTIsTAF} does not actually
depend on this paper.
In the proof of Lemma~\ref{2D}, we need to know that the
ordinary irrational rotation algebras have real rank zero,
but this follows from Remark~6 of~\cite{EE}.
In the proof of Proposition~\ref{IrratAPR}, we need to
know that traces determine order on \pj s in $A_{\te}$ whenever
$A_{\te}$ is simple.
The proof of this in~\cite{BKR} does not rely on~\cite{Hg},
and in any case an independent proof (valid whenever
$\te$ is not purely rational) is contained in~\cite{Rf1}.
And in the application of Theorem~6.4 of~\cite{Ks4} in
the proof of Proposition~\ref{SingleStep}, we use the fact
that $A_{\af}$  has a unique tracial state, rather than
real rank zero, to show that Kishimoto's conditions hold.
\end{rmk}

To finish the proof that $A_{\te}$ is an AT~algebra,
we use the following consequence of H.~Lin's classification
theorem~\cite{Ln15} for \ca s with tracial rank zero.
This result is well known,
but we have been unable to find it explicitly in the literature.
Accordingly, we give it here.
We include the AH and AF cases as well as the AT~case
for convenient reference elsewhere,
because they have the same proof,
although we do not use them here.

Recall that an AH~algebra
is a direct limit of finite direct sums of corners
of homogeneous \ca s whose primitive ideal spaces are
finite CW~complexes.
See, for example, the statement of Theorem~4.6 of~\cite{EG2},
except that we omit the restrictions there on the type of
CW~complexes which may appear;
or see~2.5 of~\cite{LnTAF}.

\begin{prp}\label{ConseqOfClass}
Let $A$ be a simple infinite dimensional
separable unital nuclear \ca\  %
with tracial rank zero and which \suct\  (Theorem~1.17 of~\cite{RSUCT}).
Then $A$ is a simple AH~algebra with real rank zero
and no dimension growth.
If $K_* (A)$ is torsion free, then $A$ is an AT~algebra.
If, in addition, $K_1 (A) = 0,$ then $A$ is an AF~algebra.
\end{prp}

\begin{proof}
Theorems~6.11 and~6.13 of~\cite{LnTTR} show that
$K_0 (A)$ is weakly unperforated
and satisfies the Riesz interpolation property
(equivalently, by Proposition~2.1 of~\cite{Gd0},
the Riesz decomposition property).
We now apply Theorem~4.20 of~\cite{EG2} to find a
unital AH~algebra $B$ with real rank zero and no dimension
growth whose ordered scaled K-theory is the same as that of $A.$
Since $A$ is simple, so is the partially ordered group $K_0 (A),$
and therefore $B$ is also simple.
If $K_* (A)$ is torsion free,
we claim that there is a simple AT~algebra $B$ with real rank zero
whose ordered scaled K-theory is the same as that of $A.$
To prove this, note that
$K_0 (A)$ can't be $\Z$ because $A$ has real rank zero;
then we apply the proof of Theorem~8.3 of~\cite{Ell2}.
(As noted in the introduction to~\cite{Ell2},
the part of the order involving $K_1$
is irrelevant in the simple case.)
We can certainly take the groups in the direct limit decomposition
to be torsion free,
so that the proof shows that all the algebras in the direct system
constructed there may be taken to have primitive ideal space
the circle or a point.
Then Theorem~4.3 of~\cite{Ell2} shows they may all be taken
to have primitive ideal space the circle.
This gives the required AT~algebra $B.$
Finally, if in addition $K_1 (A) = 0,$
following~\cite{Ef} we may find a simple AF~algebra $B$
whose ordered scaled K-theory is the same as that of $A.$

Proposition~2.6 of~\cite{LnTAF}
(with ${\mathcal{C}}$ as defined in~2.5 of~\cite{LnTAF})
implies that simple AH~algebras with real rank zero
and no dimension growth have tracial rank zero.
In particular, $B$ has tracial rank zero.
So the classification theorem for \ca s with tracial rank zero,
Theorem~5.2 of~\cite{Ln15},
implies that $A \cong B.$
\end{proof}

\begin{thm}\label{NCTIsAT}
Let $\te$ be a nondegenerate skew symmetric real $d \times d$ matrix,
with $d \geq 2.$
Then $A_{\te}$ is a simple AT~algebra with real rank zero.
\end{thm}

\begin{proof}
Using Theorem~1.17 of~\cite{RSUCT} (see the preceding discussion
for the definition of ${\mathcal{N}}$),
it follows from Lemma~\ref{ItCrPrd}
that $A_{\te}$ satisfies the Universal Coefficient Theorem.
Clearly $A_{\te}$ is separable and nuclear.
Since
\[
K_0 (A_{\bt}) \cong K_1 (A_{\bt}) \cong \Z^{2^{n - 1}}
\]
by Lemma~\ref{ItCrPrd} and by repeated application of the
Pimsner-Voiculescu exact sequence~\cite{PV},
Theorem~\ref{NCTIsTAF} and Proposition~\ref{ConseqOfClass}
imply that $A_{\te}$ is a simple AT~algebra with real rank zero.
\end{proof}

We now consider the isomorphism and Morita equivalence classification of
simple higher dimensional noncommutative tori.
For a nondegenerate skew symmetric real $d \times d$ matrix $\te,$
Elliott~\cite{El0} has determined the range of the unique tracial state
$\ta_{\te}$ acting on $K_0 (A_{\te}),$
in terms of the ``exterior exponential''
$\exp_{\wedge} (\te) \colon \Ld^{\mathrm{even}} \Z^d \to \R.$
We regard $\te$
as a linear map from $\Z^d \wedge \Z^d$ to $\R.$
Following~\cite{El0},
if $\ph \colon \Ld^k \Z^d \to \R$
and $\ps \colon \Ld^l \Z^d \to \R$
are linear, we take, by a slight abuse of notation,
$\ph \wedge \ps \colon \Ld^{k + l} \Z^d \to \R$
to be the functional obtained from the alternating functional
on $(\Z^d)^{k + l}$ defined as the antisymmetrization of
\[
(x_1, x_2, \ldots, x_{k + l})
 \mapsto \ph (x_1 \wedge x_2 \wedge \cdots \wedge x_k)
       \ps (x_{k + 1} \wedge x_{k + 2} \wedge \cdots \wedge x_{k + l}).
\]
In a similar way, we take
$\ph \oplus \ps \colon \Ld^k \Z^d \oplus \Ld^l \Z^d \to \R$
to be $(\xi, \et) \mapsto \ph (\xi) + \ps (\et).$
Then by definition
\[
\exp_{\wedge} (\te)
     = 1 \oplus \te \oplus \tfrac{1}{2} (\te \wedge \te)
            \oplus \tfrac{1}{6} (\te \wedge \te \wedge \te)
            \oplus \cdots \colon
      \Ld^{\mathrm{even}} \Z^d \to \R.
\]
Elliott's result for the nondegenerate case is then as follows.

\begin{thm}\label{T:Elliott}
Let $\te$ be a nondegenerate skew symmetric real $d \times d$ matrix.
Then there is an isomorphism
$h \colon K_0 (A_{\te_j}) \to \Ld^{\mathrm{even}} \Z^d$
such that $\exp_{\wedge} (\te) \circ h = (\ta_{\te_j})_*,$
and such that $h ([1])$
is the standard generator $1 \in \Ld^0 (\Z^d) = \Z.$
\end{thm}

\begin{proof}
This is~1.3, Theorem~2.2, and Theorem~3.1 of~\cite{El0}.
\end{proof}

For the Morita equivalence result,
we need the following lemma.

\begin{lem}\label{L:EqOfMaps}
Let $G_1$ and $G_2$ be finitely generated free abelian groups
with the same rank,
and let $f_1 \colon G_1 \to \R$ and $f_2 \colon G_2 \to \R$
be \hm s with the same range.
Then there exists an isomorphism $g \colon G_1 \to G_2$
such that $f_2 \circ g = f_1.$
\end{lem}

\begin{proof}
Let $D \S \R$ be the common range.
Then $D$ is a finitely generated subgroup of $\R,$
so is free.
Let $t_1, t_2, \ldots, t_k \in D$ form a basis.

Choose $\et_1, \et_2, \ldots, \et_k \in G_1$
such that $f_1 (\et_j) = t_j$ for $1 \leq j \leq k.$
Choose elements $\et_{k + 1}, \et_{k + 2}, \ldots, \et_n \in G_1$
which form a basis for $\Ker (f_1).$
We claim that $\et_1, \et_2, \ldots, \et_n$
form a basis for $G_1.$
To prove linear independence,
suppose $\sum_{j = 1}^n \af_j \et_j = 0.$
Apply $f_1$ to get $\af_1 = \af_2 = \cdots = \af_k = 0,$
and use linear independence
of $\et_{k + 1}, \et_{k + 2}, \ldots, \et_n.$
To show that they span $G_1,$
let $\et \in G_1,$
choose $\af_1, \af_2, \ldots, \af_k \in \Z$
such that $f_1 (\et) = \sum_{j = 1}^k \af_j t_j,$
and use $\et - \sum_{j = 1}^k \af_j \et_j \in \Ker (f_1)$
to write this element as an integer combination of
$\et_{k + 1}, \et_{k + 2}, \ldots, \et_n.$

Similarly, there is a basis
$\mu_1, \mu_2, \ldots, \mu_n$ for $G_2$
such that $g (\mu_j) = t_j$ for $1 \leq j \leq k$
and $g (\mu_j) = 0$ for $k + 1 \leq j \leq n.$
(It has the same number of elements because $G_1$ and $G_2$
have the same rank.)
The required isomorphism $g$ is now defined by
specifying $g (\et_j) = \mu_j$ for $1 \leq j \leq n.$
\end{proof}

\begin{thm}\label{MoritaClassOfNCT}
Let $\te_1$ and $\te_2$
be nondegenerate skew symmetric real $d \times d$ matrices,
with $d \geq 2.$
Then $A_{\te_1}$ is strongly Morita equivalent to $A_{\te_2}$
\ifo\  there exists $\ld > 0$ such that
$\exp_{\wedge} (\te_1)$ and $\ld \exp_{\wedge} (\te_2)$
have the same range.
\end{thm}

\begin{proof}
By Theorem~\ref{T:Elliott},
the condition is equivalent to the existence of
$\ld > 0$ such that $(\ta_{\te_1})_*$ and $\ld (\ta_{\te_2})_*$
have the same range.

The condition is certainly necessary.
For sufficiency, use Lemma~\ref{L:EqOfMaps}
to find an isomorphism $g \colon K_0 (A_{\te_1}) \to K_0 (A_{\te_2})$
such that $\ld (\ta_{\te_2})_* \circ g = (\ta_{\te_1})_*.$
Since $A_{\te_2}$ has tracial rank zero,
there are $n \in \N$ and a \pj\  $p \in M_n (A_{\te_2})$
such that $[p] = g ([1]).$
Then $A_{\te_1}$ and $p M_n (A_{\te_2}) p$ have isomorphic
Elliott invariants,
so $A_{\te_1} \cong p M_n (A_{\te_2}) p$
by Theorem~5.2 of~\cite{Ln15}.
\end{proof}

For the isomorphism classification, we have:

\begin{thm}\label{ClassOfNCT}
Let $\te_1$ and $\te_2$
be nondegenerate skew symmetric real $d \times d$ matrices,
with $d \geq 2.$
Then $A_{\te_1} \cong A_{\te_2}$ \ifo\  %
there is an isomorphism
$g \colon \Ld^{\mathrm{even}} \Z^d \to \Ld^{\mathrm{even}} \Z^d$
such that $\exp_{\wedge} (\te_2) \circ g = \exp_{\wedge} (\te_1)$
and such that $g$ sends the standard generator $1 \in \Ld^0 (\Z^d) = \Z$
to itself.
\end{thm}

\begin{proof}
By Theorem~\ref{T:Elliott}, the condition is equivalent to isomorphism
of the Elliott invariants of $A_{\te_1}$ and $A_{\te_2}.$
Apply Theorem~5.2 of~\cite{Ln15}.
\end{proof}

One might hope that it would suffice to require that
$(\ta_{\te_1})_*$ and $(\ta_{\te_2})_*$ have the same range.
We show by example that this is not true.

\begin{exa}\label{E:NonIso}
Choose $\bt, \gm \in \R$ such that $1, \bt, \gm$ are
linearly independent over $\Q.$
Set
\[
\te_1 = \left( \begin{array}{ccc}
  0     &         \bt   &  \gm               \\
- \bt   &         0     &  \frac{2}{5}       \\
- \gm   & - \frac{2}{5} &  0
\end{array} \right)
\andeqn
\te_2 = \left( \begin{array}{ccc}
  0     &         \bt   &  \gm               \\
- \bt   &         0     &  \frac{1}{5}       \\
- \gm   & - \frac{1}{5} &  0
\end{array} \right).
\]
We claim that:
\begin{enumerate}
\item\label{NonIso_Nondeg}
$\te_1$ and $\te_2$ are nondegenerate.
\item\label{NonIso_SameR}
$(\ta_{\te_1})_*$ and $(\ta_{\te_2})_*$ have the same range.
\item\label{NonIso_NonIso}
$A_{\te_1} \not\cong A_{\te_2}.$
\end{enumerate}

Set $\ld_1 = \frac{2}{5}$ and $\ld_2 = \frac{1}{5},$
giving
\[
\te_l = \left( \begin{array}{ccc}
  0     &     \bt   &  \gm         \\
- \bt   &     0     &  \ld_l       \\
- \gm   & - \ld_l   &  0
\end{array} \right).
\]

To prove~(\ref{NonIso_Nondeg}),
we verify the condition of Lemma~\ref{NondegCond}.
Thus, suppose $x \in \Q^d$ and
$\langle x, \, \te_l y \rangle \in \Q$ for all $y \in \Q^d.$
Putting $y = (1, 0, 0),$ we get $- \bt x_2 - \gm x_3 \in \Q,$
whence $x_2 = x_3 = 0.$
Putting $y = (0, 1, 0)$ and using $x_3 = 0,$
we get $\bt x_1 \in \Q,$ whence $x_1 = 0.$
So $x = 0,$ proving nondegeneracy.

For the proofs of (\ref{NonIso_SameR}) and~(\ref{NonIso_NonIso}),
apply Theorem~\ref{T:Elliott} to $\te_1$ and $\te_2,$
obtaining isomorphisms $h_1$ and $h_2.$
Use $\exp_{\wedge} (\te_l) \circ h_l = (\ta_{\te_l})_*,$
and identify $\Ld^{2} \Z^3$ with $\Z^3$ in such a way
that $\te_l,$
regarded as a linear map from $\Z^d \wedge \Z^d$ to $\R,$
sends the standard basis elements to $\ld_l,$ $\bt,$ and $\gm.$
Thus there are isomorphisms $h_l \colon K_0 (A_{\te_l}) \to \Z^4$
such that $h_l ([1]) = (1, 0, 0, 0)$
and such that the map $f_l \colon \Z^4 \to \R$
given by $n \mapsto n_1 + \ld_l n_2 + \bt n_3 + \gm n_4$
satisfies $(\ta_{\te_l})_* = f_l \circ h_l.$ 

We prove~(\ref{NonIso_SameR}) by showing that the ranges
of $f_1$ and $f_2$ are
equal to $\frac{1}{5} \Z + \bt \Z + \gm \Z.$
This is obvious for $f_2.$
Also, it is obvious that
\[
1, \tfrac{2}{5}, \bt, \gm
  \in f_1 (\Z^4)
  \subset \tfrac{1}{5} \Z + \bt \Z + \gm \Z,
\]
whence also $f_1 (\Z^4) = \frac{1}{5} \Z + \bt \Z + \gm \Z.$

We turn to the proof of~(\ref{NonIso_NonIso}).
It suffices to show that there is no isomorphism
$g \colon \Z^4 \to \Z^4$ such that $f_2 \circ g = f_1$
and $g (1, 0, 0, 0) = (1, 0, 0, 0).$
Suppose we had such a map $g.$
The equation $g (1, 0, 0, 0) = (1, 0, 0, 0)$ determines the first
column of the matrix of $g.$
Since
\[
f_2 (0, 2, 0, 0) = \tfrac{2}{5},
\,\,\,\,\,\,
f_2 (0, 0, 1, 0) = \bt,
\andeqn
f_2 (0, 0, 0, 1) = \gm,
\]
the other columns are determined by
\[
g (0, 1, 0, 0) \in (0, 2, 0, 0) + \Ker (f_2),
\,\,\,\,\,\,
g (0, 0, 1, 0) \in (0, 0, 1, 0) + \Ker (f_2),
\]
and
\[
g (0, 0, 0, 1) \in (0, 0, 0, 1) + \Ker (f_2).
\]
Now
\[
\Ker (f_2) = \{ (- r, \, 5 r, \, 0, \, 0) \in \Z^4 \colon r \in \Z \}.
\]
So there are $r, s, t \in \Z$ such that
\[
g = \left( \begin{array}{ccccccc}
  1     &  - r      &  - s   &  - t      \\
  0     &  2 + 5 r  &  5 s   &  5 t      \\
  0     &  0        &  1     &  0        \\
  0     &  0        &  0     &  1
\end{array} \right).
\]
Now $g$ can be invertible only if $2 + 5 r = \pm 1,$
which is not possible for $r \in \Z.$
This proves~(\ref{NonIso_NonIso}).
\end{exa}

We do, however, get the following result.
Recall that the opposite algebra $A^{\mathrm{op}}$ of a \ca\  $A$
is the algebra $A$ with the multiplication reversed
but all other operations, including the scalar multiplication, the same.

\begin{cor}\label{IsomOpp}
Let $\te$ be a nondegenerate skew symmetric real $d \times d$ matrix,
with $d \geq 2.$
Then $(A_{\te})^{\mathrm{op}} \cong A_{\te}.$
\end{cor}

\begin{proof}
Using classification (for example, Theorem~5.2 of~\cite{Ln15}),
one sees that
every simple AT~algebra $A$ with real rank zero is isomorphic to
its opposite algebra,
because the ordered K-theory of $A^{\mathrm{op}}$
is the same as the ordered K-theory of $A.$
\end{proof}

As far as we know,
it is unknown whether $(A_{\te})^{\mathrm{op}} \cong A_{\te}$
for general degenerate $\te.$ 

\begin{cor}\label{NoRP}
Under the hypotheses of Proposition~\ref{IrratAPR},
and with the additional condition $n \neq 1,$
the action of $\Zqn$ on $A_{\te}$ generated by $\af$
does not have the Rokhlin property.
\end{cor}

\begin{proof}
Clearly this action is trivial on $K_* (A_{\te}).$
It follows from Theorem~\ref{NCTIsAT}
that the hypotheses on the algebra in Theorem~3.5 of~\cite{Iz2}
are satisfied.
This theorem implies, in particular, that if the
action had the Rokhlin property,
then $A_{\te}$ would be isomorphic to its tensor product with
the $n^{\infty}$~UHF algebra.
The K-theory shows this is impossible.
\end{proof}

\end{document}